\documentclass[11pt]{article}
\usepackage{graphicx}
\usepackage{psfrag}
\usepackage{amsmath}
\usepackage{amssymb}
\usepackage{epsfig}
\usepackage{multicol}

\newtheorem{theorem}{Theorem}

\newtheorem{lemma}[theorem]{Lemma}

\title{Computing Boundary Slopes of 2-bridge  Links}
\author{
Jim Hoste\\
Pitzer College\\
\\
Patrick D. Shanahan\\
Loyola Marymount University
}

\begin{document}
\maketitle

\begin{abstract}
We describe an algorithm  for computing boundary slopes of 2-bridge links. As an example, we work out the slopes of the links obtained by $1/k$ surgery on one component of the Borromean rings. A table of all boundary slopes of all 2-bridge links with 10 or less crossings is also included.

\end{abstract}

\section{Introduction} In a series of papers by Hatcher and Thurston \cite{HT:1985}, Floyd and Hatcher \cite{FH:1988}, and Hatcher and Oertel~\cite{HO:1989}, the set of incompressible, boundary incompressible surfaces in the complement of a 2-bridge knot or link, or  a Montesinos knot, are completely described and classified. In the case of knots, these papers also describe the possible boundary slopes that occur and in \cite{HO:1989} a table of all boundary slopes of Montesinos knots with 10 or fewer crossings is given. Unfortunately, this table contains several errors. However a  corrected table has been published by Dunfield~\cite{D:2001}. Moreover, the computer program written and used by Dunfield is available at {\tt http://www.CompuTop.org}.

While Floyd and Hatcher explicitly describe all incompressible surfaces in the complement of a 2-bridge link,  they do not compute the boundary slopes of these surfaces,  saying only that it should be possible in principle.  In his Ph.D. thesis~\cite{L:1993}, Lash starts with their construction and shows how to compute the associated boundary slopes. His ultimate goal was to compare the set of boundary slopes of   the Whitehead link, ${\cal L}_{3/8}$, to those predicted by the algebraic-geometric machinery of Culler-Shalen \cite{CS:1983}. Applying his algorithm to ${\cal L}_{3/8}$, Lash was able to show that in this case, every boundary slope arises from a degenerating sequence of representations of the link group into $\mbox{SL}_2 \Bbb C$. Ohtsuki \cite{O:1994} has shown this to be true for all 2-bridge knots (excluding slopes that correspond to a fiber in a fibration), but the question remains open for all 2-bridge links in general. To investigate this question it would obviously be helpful to have boundary slope data for all 2-bridge links. Unfortunately, Lash's thesis has never been published and  tables of boundary slopes of links have not been available.

In this paper we describe Lash's algorithm and develop an improved algorithm that is easier to use. This allows us to describe the types of boundary slopes that can occur. As an illustration of our techniques, we compute the boundary slopes of the $(4k-1)/(8k)$ 2-bridge links. These links may also be described as  $1/k$ surgery on one component of the Borromean rings. For this class of links, we have found  an explicit description of their eigenvalue varieties and in  a forthcoming paper will investigate the relationship between the actual boundary slopes and those detected by the eigenvalue variety.

We have written a computer program to implement our algorithm and include here a table of boundary slopes of all 2-bridge links up to 10 crossings.  In our table, links through 9 crossings are also identified by their index in Rolfsen's table \cite{Ro:1990}. Our program, as well as  a much larger table to 16 crossings, will eventually be available at {\tt http://www.CompuTop.org} and as part of {\sl Knotscape}~\cite{HT:1998}. 

In Section~\ref{section:floydandhatcher} we briefly describe Floyd and Hatcher's construction for 2-bridge links. The reader is referred to their original paper for more detail. Then in Section~\ref{section:lash} we describe Lash's algorithm for finding the boundary slopes of a given 2-bridge link using Floyd and Hatcher's construction. In Section~\ref{section:improvedalgorithm} we improve the algorithm and discuss some of its theoretical consequences. The next section includes a nice example for the infinite family of 2-bridge links already mentioned above. Finally, in Section~\ref{section:table} we tabulate boundary slope data for all 2-bridge links with 10 or less crossings.

\section{Floyd and Hatcher's Construction}
\label{section:floydandhatcher}

\begin{figure}
    \begin{center}
    \leavevmode
    \scalebox{.3}{\includegraphics{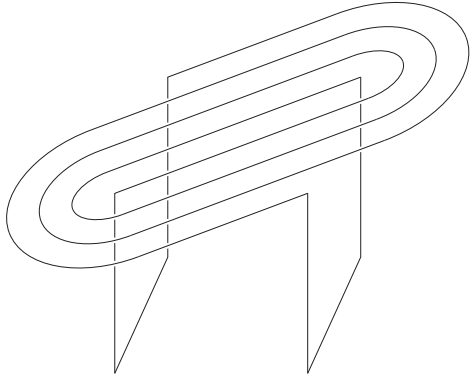}}
    \end{center}
\caption{The 2-bridge link ${\cal L}_{p/q}$, for $p=3$ and $q=4$.}
\label{figure:2-bridge link diagram}
\end{figure}

Let $p$ and $q$ be relatively prime positive integers such that $0<p<q$, $p$ is odd, and $q$ is even. We assume that the reader is familiar with the standard 2-bridge diagram of the 2-bridge link ${\cal L}_{p/q}$. For example, ${\cal L}_{3/4}$ is shown in Figure~\ref{figure:2-bridge link diagram}. It is important to note that our definition of ${\cal L}_{p/q}$ agrees with that of  \cite{HO:1989} and \cite{L:1993}, but is the mirror image of the more conventional depiction of ${\cal L}_{p/q}$ with the ``(straight) bridges on top.''  See for example, \cite{BZ:2003} or \cite{HT:1985}. 

Viewing $S^3$ as the 2-point compactification of $S^2 \times \Bbb R$, we may place ${\cal L}_{p/q}$ in $S^2 \times I$ so that it meets $S^2 \times \{0\}$ and $S^2 \times \{1\}$ each in two arcs, and each intermediate level in four points.
We may think of each level, $S^2 \times \{z\}$, as the quotient $\Bbb R^2/\Gamma$, where $\Gamma$ is the group generated by $180^\circ$ rotations of $\Bbb R^2$ about the integer lattice points $\Bbb Z^2$. The four points of the link at each intermediate level are precisely the four points of $\Bbb Z^2/\Gamma$.
The arcs at level $z=1$ are the image of  the lines in $\Bbb R^2$ which pass through integer lattice points with slope $p/q$. Similarly,   
the two arcs at level $z=0$ are the image of vertical lines through integer lattice points. Finally, PSL$_2\Bbb Z$ acts linearly on each level, leaving $\Bbb Z^2/\Gamma$ invariant.

Floyd and Hatcher next describe four basic branched surfaces, $\Sigma_A, \Sigma_B, \Sigma_C$ and $\Sigma_D$, copies of which will be stacked, one on top of the other, to build a branched surface spanning the 2-bridge link. Actually it is not these surfaces exactly, but rather homeomorphic images of them, that will be compressed vertically and stacked together. To understand the Floyd and Hatcher construction, and Lash's computation of the boundary slopes, we need to first understand these four basic building blocks.

\begin{figure}
\psfrag{a}{$\Sigma_A$}
\psfrag{b}{$\Sigma_B$}
\psfrag{c}{$\Sigma_C$}
\psfrag{d}{$\Sigma_D$}
\psfrag{x}{$\alpha$}
\psfrag{y}{$\beta$}
\psfrag{z}{$\alpha-\beta$}
\psfrag{n}{$n$}
\psfrag{q}{$\alpha-\beta-n$}
\psfrag{bn}{$\beta - n$}
\psfrag{abt}{$\frac{\alpha-\beta}{2}$}
\psfrag{r}{$0$-level}
\psfrag{s}{$\frac{1}{2}$-level}
\psfrag{t}{$1$-level}
    \begin{center}
    \leavevmode
 \includegraphics[width=5.15in]{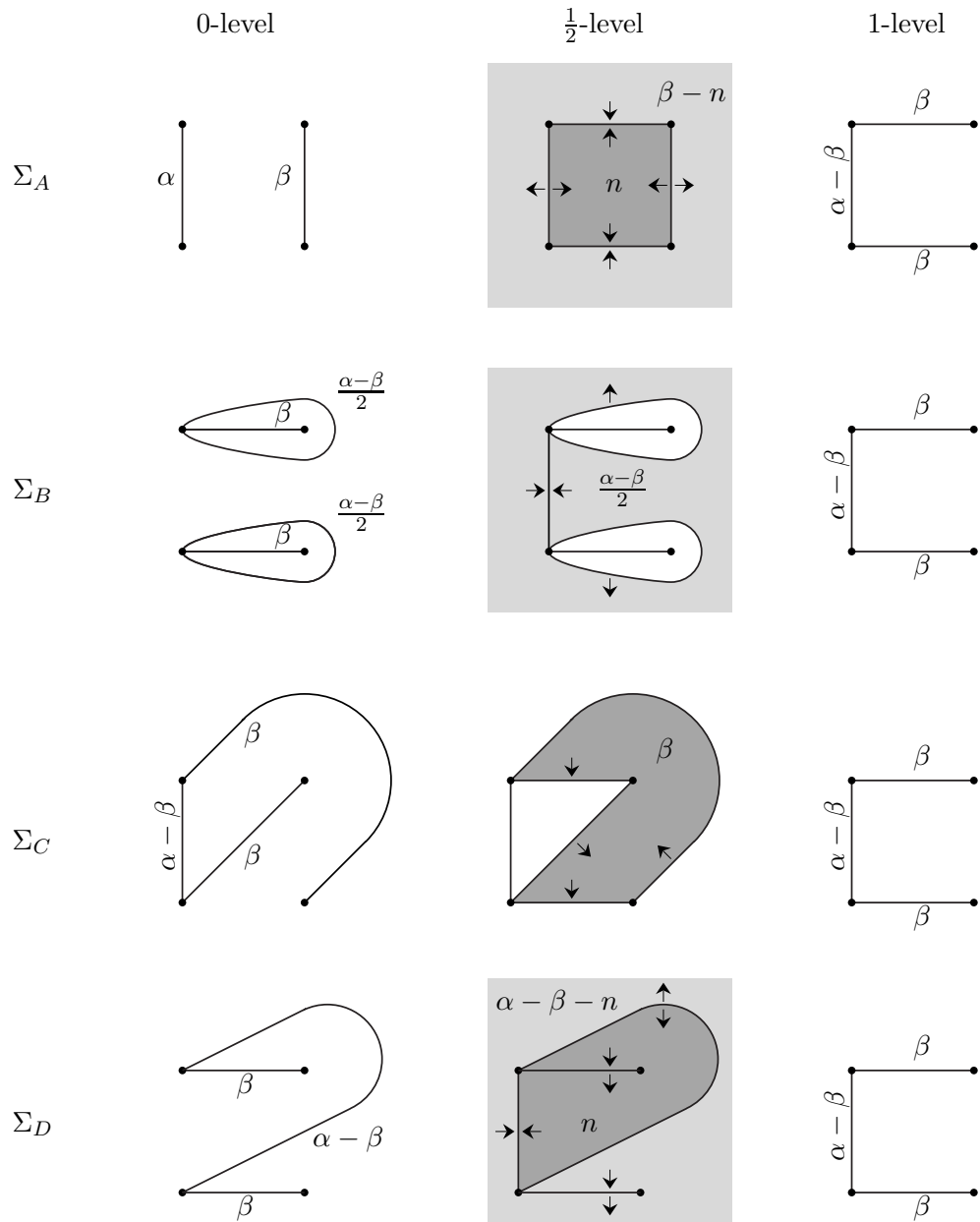}
\end{center}
\caption{The four basic branched surfaces.}
\label{figure:branchedsurfaces}
\end{figure}

Beautiful illustrations of the four surfaces are given in \cite{FH:1988} which we will not attempt to reproduce here. Instead, we describe the surfaces in a different fashion. Each  is contained in $S^2 \times I$. 
 In both of the half-intervals $[0,\frac{1}{2})$ and $(\frac{1}{2}, 1]$ the surface is a product of the half-interval with a finite number of disjoint embedded arcs in the 2-sphere with endpoints at the four points of $\Bbb Z^2/\Gamma$.  The transversality of the surface with the horizontal levels completely degenerates at the $\frac{1}{2}$-level, where branching occurs which allows the arc system at the $0$-level to transition to the arc system at the $1$-level. Figure~\ref{figure:branchedsurfaces} shows cross-sections of each branched surface at heights $0, \frac{1}{2}$, and $1$. At the $\frac{1}{2}$-level the shaded areas indicate horizontal parts of the surface. Each of the four branched surfaces can carry a variety of embedded surfaces in the usual way indicated by the number of sheets, or weights, on each piece of the branched surface. The weights are also indicated in Figure~\ref{figure:branchedsurfaces} as well as arrows at the $\frac{1}{2}$-level which indicate the direction of branching as one moves up from the $0$-level to the $1$-level. Finally, notice that all the surfaces shown in Figure~\ref{figure:branchedsurfaces} carry the implicit assumption that $\alpha > \beta$ (and sometimes that $\alpha$ and $\beta$ have the same parity).

If $g=\left (\begin{array}{cc} a &b \\ c & d \end{array}\right)$ is any element of  PSL$_2\Bbb Z$,  let $\hat g=\left (\begin{array}{cc} d&c \\ b&a \end{array}\right)$. We may take $S^2 \times I$ to itself with the homeomorphism $\hat g \times$id  and carry any one of the four basic surfaces to a new branched surface that begins and ends at arc systems  with slopes depending on $g$. Images like this of the four basic surfaces can then be joined together vertically provided they have matching arc systems where they are attached. In this way we can piece together a branched surface that begins with the two arcs of ${\cal L}_{p/q} $ at level 0 of slope $\frac{1}{0}$ and ends with the two arcs of ${\cal L}_{p/q} $ at level 1 of slope $\frac{p}{q}$.  For example, suppose we begin with the Hopf link, ${\cal L}_{1/2}$. We wish to piece together a branched surface that starts at slope $\frac{1}{0}$ and ends at slope  $\frac{1}{2}$. Starting with a copy of $\Sigma_A$ we may move from two arcs at slope $\frac{1}{0}$ to an arc system consisting of three arcs: one at slope $\frac{1}{0}$ and two at slope $\frac{0}{1}$. We may then attach to this an upside-down copy of $\Sigma_D$ which then takes us to an arc system of three arcs: two at slope $\frac{0}{1}$ and one at slope $\frac{1}{2}$. Finally, we will end with an upside-down copy of $\Sigma_A$ transformed by  $\hat g \times$id where $g=\left (\begin{array}{cc} 1&0 \\ 2&1 \end{array}\right)$. Notice that the linear transformation $\hat g$ takes lines of slope $\frac{1}{0}$ to lines of slope $\frac{1}{2}$ and  lines of slope $\frac{0}{1}$ to lines of slope $\frac{0}{1}$, because
$$  \left (\begin{array}{cc} 1 &2 \\ 0&1 \end{array}\right) \left ( \begin{array}{c} 0\\1\end{array} \right )= \left ( \begin{array}{c} 2\\1 \end{array} \right )
\quad \mbox{and} \quad
\left (\begin{array}{cc} 1 &2 \\ 0&1 \end{array}\right) \left ( \begin{array}{c} 1\\0 \end{array} \right )= \left ( \begin{array}{c}1\\ 0 \end{array} \right ).
$$
Therefore, the upside-down copy of $\Sigma_A$ transformed by $\hat g \times$id ends with the desired arc system of two arcs of slope $\frac{1}{2}$.

\begin{figure}
\psfrag{a}{$1/1$}
\psfrag{b}{$3/4$}
\psfrag{c}{$2/3$}
\psfrag{d}{$1/2$}
\psfrag{e}{$1/3$}
\psfrag{f}{$1/4$}
\psfrag{g}{$0/1$}
\psfrag{h}{$-1/2$}
\psfrag{i}{$-1/1$}
\psfrag{j}{$1/0$}
\psfrag{k}{$2/1$}
\psfrag{l}{$3/2$}

\psfrag{x}{C}
\psfrag{y}{A}
\psfrag{z}{B}
\psfrag{w}{D}

\psfrag{u}{$D_1$}
\psfrag{v}{$D_0=D_\infty$}

    \begin{center}
    \leavevmode
    \scalebox{.8}{\includegraphics{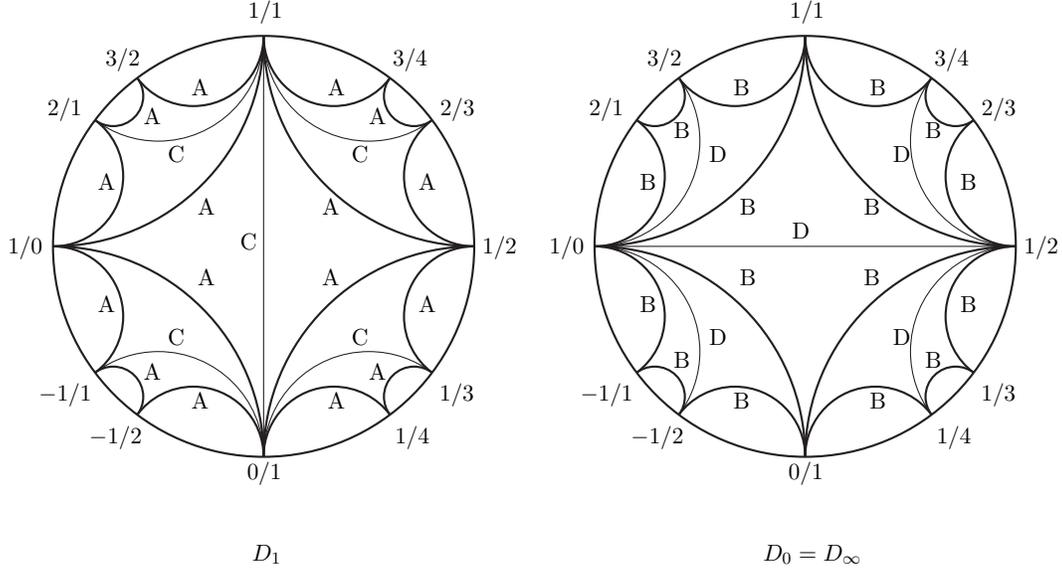}}
    \end{center}
\caption{The diagrams $D_0, D_1$, and $D_\infty$.}
\label{figure:D0andD1}
\end{figure}

There is a beautiful correspondence between branched surfaces constructed in this way and continued fraction expansions of $p/q$ which in turn may be viewed as paths in the following diagrams.
Consider first the tessellation $D_1$ of $\Bbb H^2$ by ideal triangles shown in Figure~\ref{figure:D0andD1}. The rationals, together with $\frac{1}{0}$, are arranged around the unit circle as shown, and two fractions $\frac{a}{b}$ and $\frac{c}{d}$ are connected by a geodesic if and only if $ad-bc=\pm 1$. (This diagram contains the Stern-Brocot tree generated from $\frac{1}{0}$ and $\frac{0}{1}$ by adding fractions the ``wrong way'' according to the (mis)rule $\frac{a}{b}+\frac{c}{d}=\frac{a+c}{c+d}$. See \cite{GKP:1994}, for example.) The group of orientation preserving symmetries of $D_1$ is PSL$_2\Bbb Z$. Let $G\subset \mbox{PSL}_2\Bbb Z$ be the subgroup of Mobius transformations given by $z\to \frac{az+b}{cz+d}$ where $c$ is even. It follows that the  triangle $\{\frac{1}{0}, \frac{0}{1}, \frac{1}{1}\}$ is a fundamental domain for the action of $G$ and the $G$-images of the ideal quadralateral $Q=\{\frac{1}{0}, \frac{0}{1},\frac{1}{2}, \frac{1}{1}\}$ tessellate $\Bbb H^2$. If we delete the $G$-orbit of the diagonal $\{\frac{0}{1}, \frac{1}{1}\}$ and replace it with the $G$-orbit of the opposite diagonal $\{\frac{1}{0}, \frac{1}{2}\}$, we obtain a new diagram called $D_0$, which is also shown in Figure~\ref{figure:D0andD1}.

\begin{figure}
\psfrag{a}{$1/1$}
\psfrag{b}{$1/2$}
\psfrag{c}{$0/1$}
\psfrag{d}{$1/0$}

\psfrag{e}{A}   
\psfrag{f}{B}   
\psfrag{g}{C}   
\psfrag{h}{D} 
\psfrag{e0}{A$_0$}   
\psfrag{f0}{B$_0$}   
\psfrag{g0}{C$_0$}   
\psfrag{h0}{D$_0$}
\psfrag{r}{D$_1$}
\psfrag{s}{D$_t$}
\psfrag{t}{D$_0$}   
  
\begin{center}
    \leavevmode
    \scalebox{.85}{\includegraphics{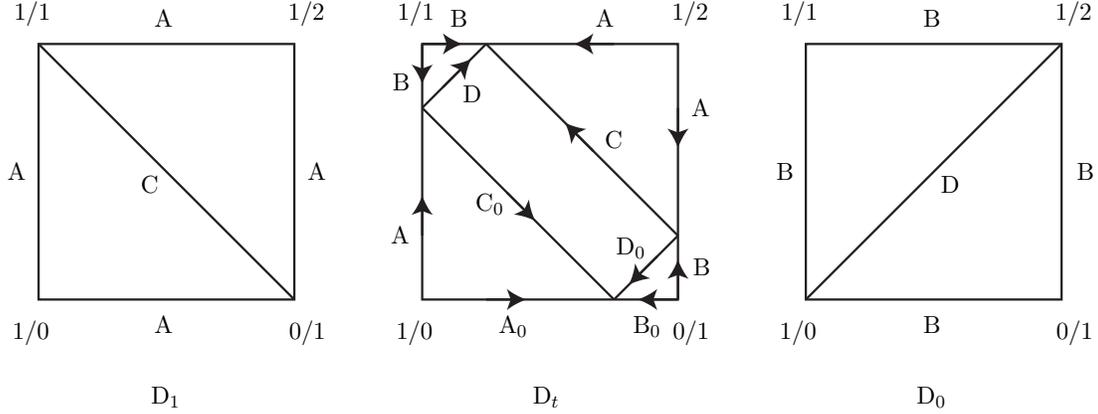}}
    \end{center}
\caption{Expanding the diagonals of $D_0$ and $D_1$ to obtain $D_t$
\
(pictured here with $t=3/4$).}
\label{figure:quadInDt}
\end{figure}

Between $D_0$ and $D_1$ there exists a family of oriented diagrams $D_t$, for $0<t<1$, obtained by expanding each diagonal in $D_1$ (labeled C in Figure~\ref{figure:D0andD1}) to a rectangle which may then be collapsed to the opposite diagonal (labeled D in  Figure~\ref{figure:D0andD1}), thus giving $D_0$. The intersection of each of $D_0, D_t$, and $D_1$ with the fundamental quadralateral $Q$ is shown in Figure~\ref{figure:quadInDt}.  The edges of $D_t$ fall into four $G$-orbits which are named $A, B, C$, and $D$ and which have, respectively, representative edges $A_0, B_0, C_0$ and $D_0$ as defined in Figure~\ref{figure:quadInDt}. As $t\to 0$ the edges degenerate into $B$ and $D$-type edges in $D_0$. If instead, $t\to 1$, the edges degenerate into $A$ and $C$-type edges in $D_1$. We may orient the edges of $D_t$ by using the orientations of $A_0, B_0, C_0$ and $D_0$  shown in Figure~\ref{figure:quadInDt}, but there is no coherent way to orient the edges of $D_0$ or $D_1$. Finally, by setting $D_t=D_{1/t}$ for $1\le t \le \infty$, we obtain a diagram for every $t\in[0,\infty]$.

Floyd and Hatcher show that the diagram $D_t$ provides a beautiful way of describing all the incompressible surfaces in the 2-bridge link exteriors. In particular, {\it minimal edge paths} in $D_t$ from $\frac{1}{0}$ to $\frac{p}{q}$ will correspond to branched surfaces which in turn will carry the incompressible surfaces. An edge  path in $D_t$ is {\it minimal} if it never contains two consecutive edges that lie in the same triangle or rectangle of $D_t$. It is not hard to see that a minimal edge path in $D_t$ (with $t \not\in \{0,1, \infty\}$) will collapse to a minimal edge path in either $D_0$ or $D_1$ as $t$ approaches 0 or 1 respectively. Moreover, these limiting paths in $D_0$ and $D_1$ uniquely determine the path in $D_t$. For a particular fraction $\frac{p}{q}$ there can only be a finite number of minimal edge paths connecting it to $\frac{1}{0}$. This follows from the fact that these  minimal paths are all contained in a unique minimal chain of quadralaterals consisting of  $Q$ and a finite number of its translates under $G$.

Each minimal edge path in $D_t$ (with $t \not\in \{0,1,\infty\}$) from $\frac{1}{0}$ to $\frac{p}{q}$ provides a recipe for piecing together images of the four basic surfaces $\Sigma_A, \Sigma_B, \Sigma_C$, and $\Sigma_D$ as follows. Suppose $\gamma$ is a path consisting of consecutive edges $e_1, e_2, \dots, e_m$. For each edge $e_i$ there exists an element $g_i \in G$ taking $E_i$ to $e_i$, where  $E_i$ is  one of the four representative edges $A_0, B_0, C_0$, or $D_0$ in $Q$.  Let $S_i$ be $\Sigma_A, \Sigma_B, \Sigma_C$, or $\Sigma_D$  depending on whether $E_i$ is  $A_0, B_0, C_0$, or $D_0$ respectively. Now take $S_i$ to $(\hat g_i \times \mbox{id})(S_i)$, rescaling vertically so as to place the image in $S^2 \times [\frac{i-1}{m}, \frac{i}{m}]$. Moreover, if the orientation of $g_i(E_i)$ is opposite that of $e_i$, we first reflect $S_i$ through the 2-sphere $S^2 \times \{\frac{1}{2}\}$ before applying $\hat g_i \times \mbox{id}$.

If $t$ is rational and equal to the reduced fraction $\alpha/\beta$, we can further use this information to assign weights to the branched surface as follows. If $t>1$ the branched surface is weighted as shown in Figure~\ref{figure:branchedsurfaces}. If instead, $0<t<1$  the above construction is altered by first rotating each of the basic surfaces $\Sigma_A, \Sigma_B, \Sigma_C$, and $\Sigma_D$ through $180^\circ$ in the obvious way so as to interchange the two components of ${\cal L}_{p/q}$. We then swap $\alpha$ and $\beta$ and proceed as before. Thus for every minimal edge path in $D_t$ with $t$ a positive rational number different from $1$, we have associated a weighted branched surface and thus an actual surface in the complement of the link.

If $t=0, 1$, or $\infty$ we must modify this recipe slightly. If we let $t$ approach zero or infinity, then the above constructions will limit at surfaces for which $\alpha=0$ or $\beta=0$ respectively.  But the limiting surfaces which arise from the above constructions when $\alpha \to \beta$ do not give all desired surfaces with $\alpha=\beta$. Instead, the process must be modified slightly. As $t \to 1$ the minimal edge paths in $D_t$ approach minimal edge paths in $D_1$ that consist entirely of $A$ or $C$-type edges. Returning to Figure~\ref{figure:branchedsurfaces}, notice that if $\alpha-\beta=0$, then $\Sigma_A$ with $n=\alpha$ and $\Sigma_C$ are isotopic. Thus  if $C$-type edges are involved, then $\Sigma_A$ allows for more general branching than $\Sigma_C$ and for this reason we replace each use of $\Sigma_C$ with $\Sigma_A$.

Finally, the main result of Floyd and Hatcher is the following theorem.

\begin{theorem}[Floyd and Hatcher] The orientable incompressible and meridionally incompressible surfaces in $S^3-{\cal L}_{p/q}$, without peripheral components, are exactly (up to isotopy) the orientable surfaces carried by the collection of branched surfaces associated to minimal edge paths in $D_t$ from   $\frac{1}{0}$ to $\frac{p}{q}$ and with $t \in [0, \infty]$.
\end{theorem}

Recall from \cite{FH:1988} that a surface $S$ in the complement of the link $L$ is {\it meridionally incompressible} if whenever there is a disk $D\subset S^3$ with $D\cap S=\partial D$ and such that $D$ meets $L$ transversely in one interior point, then there exists a disk $D^\prime\subset S \cup L$ with $\partial D= \partial D^\prime$ and such that $D^\prime$ also meets $L$ transversely in a single interior point.

Floyd and Hatcher then go on to explicitly describe when two surfaces constructed in this way are isotopic. However, they do not compute the boundary slopes of these surfaces,  saying only that it should be possible in principle. 

\section{Lash's Algorithm}
\label{section:lash}

The boundary of a branched surface derived from the Floyd-Hatcher construction  defines a train track on the boundary of the regular neighborhood of the link. Thus the boundary of any incompressible surface carried by the branched surface is carried by this train track. Lash's first step is to determine the train tracks for each of the four basic surfaces $\Sigma_A, \Sigma_B, \Sigma_C$ and $\Sigma_D$.

Before doing this,  we introduce some notation. Let the four points of $\Bbb Z^2/\Gamma$ be $(0,0), (1,0), (0,1)$ and $(1,1)$. At the 0-level, assume that the two arcs of ${\cal L}_{p/q}$ join the points $(0,0)$ to $(0,1)$ and $(0,1)$ to $(1,1)$. Furthermore, let ${\cal L}_{p/q}=\{K_1, K_2\}$ where $K_1$ contains $(0,0)$ and $(0,1)$. Orient $K_1$ so that it runs vertically  upward from $(0,0)$ and orient $K_2$ so that it runs vertically upward from $(1,1)$. Choose as a fundamental domain of $\Bbb R^2/\Gamma$ the region ${\cal D}=[0,1] \times [-1/2, 3/2]$. Removing small disks of radius $\epsilon$ centered at the four points of $\Bbb Z/\Gamma$ will remove semi-disks from the fundamental domain $\cal D$.  These correspond to meridional cross sections of
the regular neighborhood of the link. 

Let $\mu_1$ be the oriented meridian of $K_1$ having linking number $+1$ with $K_1$. Let $\lambda_1$ be the oriented  longitude of $K_1$ defined as follows. Start with  the line segment from $(0,1-\epsilon, \epsilon)$ to  $(0,\epsilon, \epsilon)$. Join to this the vertical segments $\{(0,\epsilon)\} \times [ \epsilon, 1-\epsilon]$ and $\{(0,1-\epsilon)\} \times [ \epsilon, 1-\epsilon]$.  Next, add the curve in $\Bbb R^2 \times \{1-\epsilon\}$ which starts at $(0,\epsilon, 1-\epsilon)$, 
ends at $(0,1-\epsilon, 1-\epsilon)$, and is parallel to $K_1$. Finally, orient $\lambda_1$ parallel to $K_1$. 

The $180^\circ$ rotation of Figure~\ref{figure:2-bridge link diagram} about the vertical axis $\{(1/2, 1/2)\} \times \Bbb R$ interchanges the components of the link, and  preserves their orientations.  
We  define the oriented meridian $\mu_2$ and longitude $\lambda_2$ of $K_2$ as the images of  $\mu_1$ and $\lambda_1$ respectively under this rotation.

We will initially compute all boundary slopes with respect to the basis $\{\mu_i, \lambda_i \}$. However, $\lambda_i$ is not necessarily a preferred longitude of $K_i$ so it will be necessary later to know its linking number with $K_i$. It is a straightforward exercise to compute this. We obtain
\begin{equation}
\mbox{lk}(K_i, \lambda_i)=-\sum_{j=1}^{\frac{q-2}{2}}(-1)^{[\![2 jp/q]\!]}
\label{eqn:linkingnumber}
\end{equation}
where $[\![x]\!]$ is the greatest integer less than or equal to $x$. Note, for example,  that $\lambda_i$ is never the preferred longitude if $q$ is a multiple of 4 since the linking number must be odd in this case.

\begin{figure}
\psfrag{00}{$(0,0)$}
\psfrag{01}{$(0,1)$}
\psfrag{10}{$(1,0)$}
\psfrag{11}{$(1,1)$}
\psfrag{a}{\small$\alpha$}
\psfrag{b}{\small$\beta$}
\psfrag{n}{\scriptsize$n$}
\psfrag{bn}{\scriptsize$\beta-n$}
\psfrag{i}{$\infty$}
\psfrag{ni}{$-\infty$}
\psfrag{ab}{\small$\alpha-\beta$}
\psfrag{ab2}{\scriptsize$\displaystyle \frac{\alpha-\beta}{2}$}
\psfrag{0}{$0$}
\psfrag{.5}{$.5$}
\psfrag{1}{$1$}
\psfrag{SA}{$\Sigma_A$}
\psfrag{SB}{$\Sigma_B$}
\psfrag{SC}{$\Sigma_C$}
\psfrag{SD}{$\Sigma_D$}
\psfrag{abn}{\scriptsize$\alpha-\beta-n$}
    \begin{center}
    \leavevmode
    \scalebox{.8}{\includegraphics{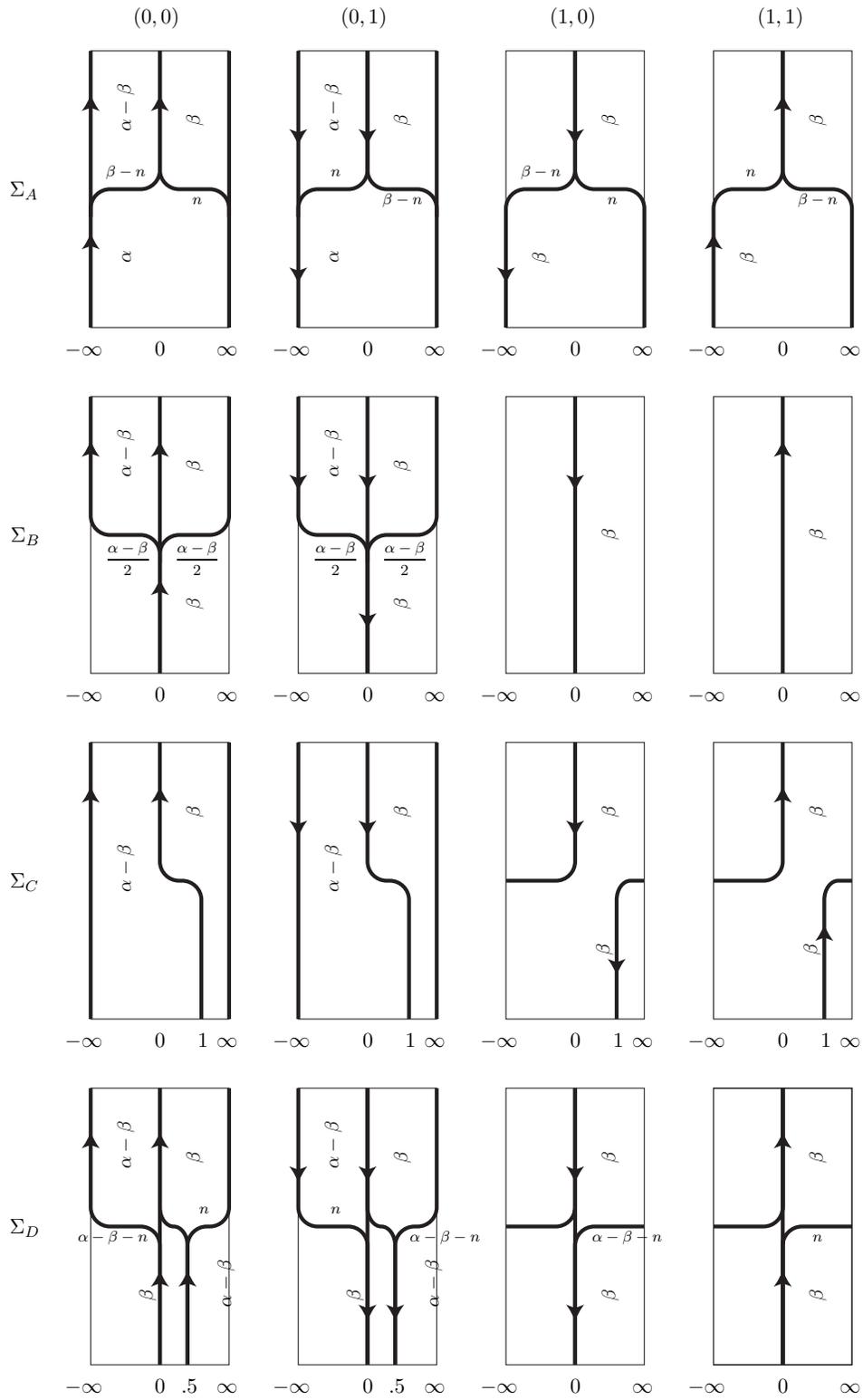}}
    \end{center}
\caption{Train track boundaries of each of the four branched surfaces.}
\label{fig:traintracks}
\end{figure}

Figure~\ref{fig:traintracks} depicts the train track boundaries of each of the four basic branched surfaces. Each row of the figure shows the train tracks on the boundary of each of the four ``columns'' which are the regular neighborhoods of the vertical segments of the link. These are depicted in the corresponding regions of ${\cal D} \times [0,1]$, that is, the product of the semicircular arcs surrounding each integer lattice point with the unit interval $[0,1]$. Notice that we have used slopes ranging from $-\infty$ to $\infty$ to parameterize the semicircular arcs.  Thus we see the train tracks for $\Sigma_A$, on each of the four columns, begin at slopes of $\pm \infty$ and end at slopes of $\pm \infty$ and $0$. Similarly, the train tracks for $\Sigma_C$, on all four columns, have curves that begin at slopes of  $1$ and end at slopes of $0$.

Suppose $\Sigma$ is a branched surface obtained by piecing together homeomorphic copies of the four basic branched surfaces and that $S$ is a surface carried by $\Sigma$. Because $\Sigma$ corresponds to a path in $D_t$ from $\frac{1}{0}$ to $\frac{p}{q}$, $\Sigma$ will always begin with a copy of $\Sigma_A$ and end with an upside-down copy of $\Sigma_A$. Let the initial weights on $\Sigma$ at the 0-level be $\alpha$ and $\beta$ as shown in Figure~\ref{figure:branchedsurfaces}. Thus, at the 0-level on the neighborhood of $K_1$, the boundary of $S$ consists of $\alpha $ arcs. If we orient one of these arcs in the direction of $K_1$ and follow it upward, it will follow the train track up  the $(0,0)$ column until it reaches the top, traverse a curve of slope $\frac{p}{q}$ at the 1-level, and then follow the train track down the $(0,1)$ column. In Figure~\ref{fig:traintracks} we have chosen to orient the train tracks parallel to $K_1$ and $K_2$ for this reason. When we return to our starting point at the bottom we may end there, or perhaps continue to travel around again in the longitudinal direction. If the boundary of $S$ consists of several components, we may orient them all parallel to $K_1$ and $K_2$  (even if these orientations are not compatible with an orientation of $S$).

Suppose $\partial S$ has $k$ components  on the regular neighborhood of $K_1$ and that each has algebraic (and geometric) intersection of $l_i$ with $\mu_1$. Thus $$l_1+l_2+\dots+ l_k=\alpha.$$ Furthermore, suppose each has algebraic intersection $m_i$ with $\lambda_1$. Then $\gcd(m_i, l_i)=1$ for all $i$ and moreover, because the components of $\partial S$ are all disjoint and nontrivial, each boundary slope $m_i/l_i$ is the same. From this it follows that $l_i$ and $m_i$ are constant, say $l_i=l$ and   $m_i=m$. Thus to determine the boundary slope $m/l$ we need only compute the total algebraic intersection $M_1=m_1+m_2 + \dots + m_k$ of $\partial S$ with $\lambda_1$ and divide by $\alpha$. Similarly, the boundary slope of $S$ on the regular neighborhood of $K_2$ is $M_2/\beta$ where $M_2$ is the total algebraic intersection of $\partial S$ with $\lambda_2$. Again, these boundary slopes are with respect to the basis $\{\mu_i, \lambda_i\}$.

To compute $M_1$ we must sum the algebraic intersection of the oriented train tracks with $\lambda_1$ on columns $(0,0)$ and $(0,1)$. We may do this one section at a time, with each section coming from one of the four basic surfaces. We may simplify matters by pulling back $\lambda_1$ under the inverse of $\hat g \times$id and counting its intersection with the standard train tracks shown in Figure~\ref{fig:traintracks}, rather than examining the images of the standard train tracks under $\hat g \times$id. If  $g=\left (\begin{array}{cc} a &b \\ c & d \end{array}\right)$ then  $\hat g^{-1}=\left (\begin{array}{cc} a & -c \\ -b & d \end{array}\right)$ and the vector
$\left( \begin{array}{c} 0 \\ 1 \end{array} \right)$, 
which represents a slope of $\infty$ (and the longitude $\lambda_1$), pulls back to the vector $\left( \begin{array}{c} -c \\ d \end{array} \right)$, which represents a slope of $-d/c$. Thus we must examine how the standard train tracks intersect the vertical line located at slope $-d/c$.

For example, suppose $\gamma$ is a minimal edge path in $D_t$ beginning at $\frac{1}{0}$ and ending at $\frac{p}{q}$. Let $e_i$ be an edge of $\gamma$ that is the image of $E_i=A_0$ under the element $g_i$ of $G$. Furthermore, suppose that the orientations of $e_i$ and $g_i(A_0)$ agree.  If $0<-d/c<\infty$ then the algebraic intersection on column $(0,0)$ is $-n$ and on column $(0,1)$ is $-(\beta-n)$. Note that on column $(0,0)$ the meridian runs from $0$ towards $\infty$ while on column $(0,1)$ the meridian runs from $0$ to $-\infty$. This explains the minus signs in the above calculations. The total contribution to $M_1$ at this level, in this case, is thus $-\beta$. If instead, $-\infty<-d/c<0$, we obtain $\beta-n$ on column $(0,0)$ and $n$ on column $(0,1)$ for a total contribution to $M_1$ of $\beta$. Note that $-d/c$ cannot be zero. This is because $c$ is even and therefore both $a$ and $d$ must be odd since $\det(g)=1$. However, $-d/c$ may equal $\pm \infty$ if $c=0$. In this case the train track is not transverse to the (pullback of the) longitude. When this occurs, we may istotope the longitude slightly by pushing it in the positive direction of $\mu_1$. So in this case, if $-d/c=\pm \infty$ we obtain an intersection on column $(0,0)$ of $\beta -n$ and on column $(0,1)$ of $-(\beta-n)$ for a net contribution to $M_1$ of zero. Table~\ref{table:m1andm2} lists these results together with the contributions to $M_1$ for the other three types of surfaces, all in the case where the orientations of $e_i$ and $g_i(E_i)$ agree. If these orientations are opposite, then we must flip each surface upside-down and it is easy to see that the effect is to negate the entries in the table.

A similar examination of the standard train tracks allows us to compute $M_2$. However in this case an additional consideration seems necessary. For any $g \in G$, $\hat g$ takes $(0,0)$ to itself and $(0,1)$ to itself since $c$ is even and both $a$ and $d$ are odd. But if $b$ is odd, then  $(1,0)$ and $(1,1)$ will be traded while if $b$ is even these lattice points will each be fixed. Thus it seems necessary to consider these cases separately. However, a close examination of the case when $b$ is odd reveals that both the orientations of the train tracks as well as the orientation of $\mu_2$ are reversed and the total contribution to $M_2$ is unaffected. Thus separate formulae for $b$ even or odd are not needed.  The results for each of the four basic surfaces are listed in Table~\ref{table:m1andm2}.

\begin{table}[htdp]
\renewcommand{\arraystretch}{1.2}
\caption{Contributions to $M_1$ and $M_2$ according to surface type.}
\label{table:m1andm2}
\begin{center}
\begin{tabular}{|c|c|c|c|}

\hline
Surface Type & $M_1$ &  $M_2$&$-d/c$\\
\hline
\hline
$A$&
\begin{tabular}{c} 
$\beta$\\ 
$-\beta$\\
$0$\\
\end{tabular} &
\begin{tabular}{c}
$\beta$\\
$-\beta$\\
$0$\\
\end{tabular} &
\begin{tabular}{c}
$-\infty<-d/c<0$\\
$0<-d/c<\infty$\\
$-d/c=\pm\infty$\\
\end{tabular} \\
 
\hline
$B$&
\begin{tabular}{c} 
$\beta-\alpha$\\ 
$\alpha-\beta$\\
$0$\\
\end{tabular} & 
\begin{tabular}{c} 
$0$\\ 
$0$\\ 
$0$\\ 
\end{tabular} &
\begin{tabular}{c} 
$-\infty<-d/c<0$\\ 
$0<-d/c<\infty$\\
$-d/c=\pm\infty$\\
\end{tabular}  \\
\hline

$C$&

\begin{tabular}{c} 
$-2\beta$\\ 
$0$\\
\end{tabular} &
\begin{tabular}{c} 
$0$\\ 
$2 \beta$\\
\end{tabular} &
\begin{tabular}{c}
$0<-d/c<1$\\ 
otherwise\\
\end{tabular}\\
\hline

$D$&
\begin{tabular}{c} 
$\beta-\alpha$\\ 
$0$\\
$\alpha-\beta$\\
\end{tabular} &
\begin{tabular}{c} 
$\alpha-\beta$\\
$\alpha-\beta$\\
$\alpha-\beta$\\
\end{tabular} &
\begin{tabular}{c} 
$-\infty<-d/c<1/2$\\ 
$-d/c=1/2, \pm\infty$\\
$1/2<-d/c<\infty$\\
\end{tabular} \\
\hline

$A$ for $C$&
\begin{tabular}{c} 
$2(\beta-n)$\\ 
$-2n$\\
\end{tabular} 
&
\begin{tabular}{c}
$2n$\\
$2(n-\beta)$\\
\end{tabular} 
&
\begin{tabular}{c}
$-\infty<-d/c<0$\\
$0<-d/c<\infty$\\
\end{tabular}\\
\hline

\end{tabular}
\end{center}
\end{table}

The data in the first four rows of Table~\ref{table:m1andm2} are sufficient to compute the boundary slopes when $t=\alpha/\beta \not\in \{0,1,\infty\}$ and even in the limiting cases $t\to 0$ or $t \to \infty$. But as mentioned already, when $t \to 1$,  if any  $C$-type edges appear in the limiting minimal path in $D_1$, then  we  replace the corresponding $C$-type surfaces with the more general $A$-type surfaces. Suppose $e$ is such a $C$-type edge in $D_1$. Then $e$ joins $a/c$ and $b/d$ in $D_1$ where both $c$ and $d$ are odd, $ad-bc=1$, and thus, exactly one of $a$ or $b$ is even. Therefore we may choose $g=\left (\begin{array}{cc} a &b \\ c & d \end{array}\right)$ and assume that $b$ is even, while $a,c$, and $d$ are all odd. Now $\hat g$ fixes $(0,0)$ and $(1,0)$ while trading $(0,1)$ and $(1,1)$. To properly use  Figure~\ref{fig:traintracks} now, we must orient columns $(0,0)$ and $(0,1)$ up and the other two down. Furthermore, the meridian $\mu_1$ points from $0$ to $\infty$ in the $(0,0)$ figure and from $0$ to $-\infty$ in the $(1,1)$ figure. Similarly, the meridian $\mu_2$ points from $0$ to $\infty$ in the $(0,1)$ figure and oppositely in the $(1,0)$ figure. It is now a simple matter to compute the contributions to $M_1$ and $M_2$ given in the last row of Table~\ref{table:m1andm2}. Here we see that the number $n$ of horizontal sheets does not cancel from the calculations. Notice also that because both $c$ and $d$ are odd, $-d/c$ can not equal zero or infinity.

Everything is now in place to compute the boundary slopes for a given 2-bridge link, ${\cal L}_{p/q}$. Beginning with $0<\beta<\alpha$, and hence $1<t<\infty$, we first find all minimal paths in $D_t$ from $1/0$ to $p/q$. For each of these paths, each edge $e$ must be identified as the image of $E \in \{A_0, B_0, C_0, D_0\}$ by some element $g=\left (\begin{array}{cc} a &b \\ c & d \end{array}\right) \in G$. Next, using $-d/c$ and the data in Table~\ref{table:m1andm2}, the contributions to $M_1$ and $M_2$ are computed and these are added over all edges in the path to obtain $M_1$ and $M_2$. Note that the entries in the table must be negated if the orientation of $e$ does not match that of  $g(E)$.  Notice also that $M_1$ and $M_2$ are always integer linear combinations of $\alpha$ and $\beta$. The boundary slope of this surface is then 
$$\frac{M_1}{\alpha}=\frac{x \alpha +y \beta}{\alpha}=x+\frac{y}{t}$$ 
on the boundary associated to $K_1$ and 
$$\frac{M_2}{\beta}=\frac{z \alpha +w \beta}{\beta}=z t+w$$
 on the boundary associated to $K_2$. This gives a 1-parameter family of boundary slopes for each rational number $t$ greater than 1. As $t$ approaches $\infty$ we may simply take the limits of the slopes obtained so far provided $z\ne0$. If $z=0$, then $M_2$ and $\beta$ both approach zero as $t$ approaches $\infty$. This means that the surface has no intersection with the regular neighborhood of $K_2$ and thus has no boundary slope associated with $K_2$.  At the other extreme, letting $t$ approach 1 will produce legitimate boundary slopes, but not all possible slopes with $\alpha=\beta$. Instead we must consider the limiting minimal edge path in $D_1$, swap $A$ for $C$-type edges (if there are any) and now use the data from the last row of Table~\ref{table:m1andm2}. Finally,  to allow for $t<1$, we must rotate all our surfaces 180$^\circ$ around the axis $\{(1/2, 1/2)\} \times \Bbb R$, thus trading $\alpha$ with $\beta$. The 4-tuple $(M_1(\alpha, \beta), \alpha, M_2(\alpha, \beta), \beta)$ of the algebraic intersections with $\lambda_1, \mu_1, \lambda_2$, and $\mu_2$ respectively is then transformed to  $(M_2(\beta, \alpha), \alpha, M_1(\beta, \alpha), \beta)$ with it now the case that $\alpha < \beta$.  Letting $t$ approach zero now corresponds to letting $\alpha$ approach zero. Finally, all of these computations are with respect to the bases  $\{\mu_1, \lambda_1\}$ and $\{\mu_2, \lambda_2\}$. To convert to a preferred bases, we must consider the linking number $l=$lk$(\lambda_1, K_1)=$lk$(\lambda_2, K_2)$ given in Equation~\ref{eqn:linkingnumber}. Converting to the preferred basis sends the 4-tuple $(M_1, \alpha, M_2, \beta)$ to $(M_1+l \alpha, \alpha, M_2+l \beta, \beta)$.

\section{An Improved Algorithm}
\label{section:improvedalgorithm}

As mentioned already, $M_1$ and $M_2$ are always integer linear combinations of $\alpha$ and $\beta$. After implementing Lash's algorithm on a computer and looking at sample data, the conclusions of the following theorem were apparent. In searching for a proof, we were led to a revision of Lash's algorithm that is much simpler to apply by hand and implement on a computer. The revised approach will be described in the proof and illustrated in the next section with an interesting example.

Before stating some results we define $M(\gamma)$ to be the pair of intersection numbers $(M_1, M_2)$ associated to the path $\gamma$ in $D_t$.

\begin{theorem}
\label{thm:xyz}
Given any path $\gamma$ (not necessarily minimal) in $D_t$ with $1<t<\infty$,   which begins and ends at vertices  of $D_1$,  $M(\gamma)$ is of the form 
$$M(\gamma)=(x \alpha+y\beta, y\alpha+z\beta)$$
where  $x \equiv z \mod 2$. If additionally, $\gamma$ begins at $\frac{1}{0}$ and ends at $\frac{p}{q}$, then $x+y \equiv 1+ q \mod 2$.

\end{theorem}

The proof of this theorem is an easy consequence of the following two lemmas.

\begin{lemma} 
\label{lemma:AandBonly}
Suppose $\gamma_t$ is a path in $D_t$ starting at $\frac{p_0}{q_0}$, ending at $\frac{p_n}{q_n}$ and consisting solely of $A$ and $B$-type edges. Letting $t \to 1$, $\gamma_t$ will collapse to the path 
$$\gamma_1=\left \{\frac{p_0}{q_0}, \frac{p_1}{q_1}, \dots, \frac{p_n}{q_n} \right \}$$
in $D_1$ containing only $A$-type edges. Then
$$\begin{array}{rl}
M(\gamma_t)&=\left ( \sum_{i=0}^{n-1}\delta_i \right)(\alpha, \beta) \mbox{ and }\\ 
M(\gamma_1)&=\left ( \sum_{i=0}^{n-1}\delta_i \right)(\beta, \beta)
\end{array}
$$
where
$$\delta_i=\left \{ \begin{array}{ll} 0, & \mbox{if  } q_iq_{i+1}= 0;\\
 p_i q_{i+1}-p_{i+1}q_i, & \mbox{otherwise.}  \end{array}\right.$$
\end{lemma}
{\bf Proof:} Let $e$ be an edge of $\gamma_1$ oriented from $\frac{a}{c}$ to $\frac{b}{d}$ where $a, b, c$ and $d$ are all positive, $c$ is even and $ad-bc=\pm1$. The corresponding matrix $g \in G$ is either  $g=\left (\begin{array}{cc}  a&b\\c&d\\ \end{array}\right  ) $ or $g=\left (\begin{array}{cc}  a&-b\\c&-d\\ \end{array}\right  ) $, whichever one has determinant one. In $\gamma_t$, $e$ begins with an $A$-type edge, oriented forward, and ends with a $B$-type edge, oriented backwards. Referring to Table~\ref{table:m1andm2}, we see that the contribution to $(M_1, M_2)$ is either $(0,0)$ if $c=0$, or $(ad-bc)(\alpha, \beta)$. It is easy to check  the remaining case where $c$ is odd and $d$ is even. 
\hfill $\Box$

\begin{figure}
\psfrag{a}{$\frac{a+b}{c+d}$}
\psfrag{b}{$\frac{a+2 b}{c+2 d}$}
\psfrag{c}{$\frac{b}{d}$}
\psfrag{d}{$\frac{a}{c}$}

\psfrag{e}{A}   
\psfrag{f}{B}   
\psfrag{g}{C}   
\psfrag{h}{D} 
\psfrag{e0}{A}   
\psfrag{f0}{B}   
\psfrag{g0}{C}   
\psfrag{h0}{D}
\psfrag{s}{D$_t$}

\psfrag{r0}{$R_0$}
\psfrag{r1}{$R_1$}
\psfrag{r2}{$R_2$}
\psfrag{r3}{$R_3$}
\psfrag{r4}{$R_4$}
  
\begin{center}
    \leavevmode
    \scalebox{.85}{\includegraphics{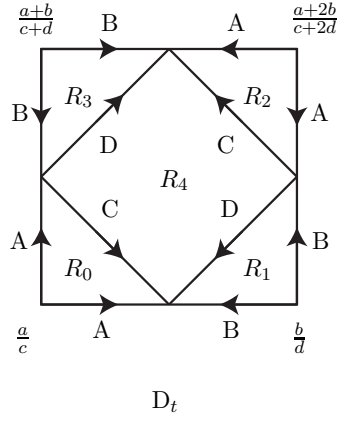}}
    \end{center}
\caption{Lemma~\ref{lemma:pathmovest} describes the result of travelling around each region.}
\label{fig:gQ}
\end{figure}

\begin{lemma}
\label{lemma:pathmovest}
 Let  $g=\left (\begin{array}{cc}  a&b\\c&d\\ \end{array}\right  ) $ be any element of $G$ and label the regions of $g(Q)$ as shown in Figure~\ref{fig:gQ}. If $\gamma_i=\partial R_i$, oriented counterclockwise, then $M(\gamma_i)$  are given by:
$$\begin{array}{l|l}
\gamma_i&M(\gamma_i)\\
\hline
\gamma_0 \mbox{ or } \gamma_2& (0, -2\beta)\\
\gamma_1 \mbox{ or } \gamma_3& (-\alpha+\beta, \alpha-\beta)\\
\gamma_4 & (-2\beta, -2\alpha+4\beta)\\
\end{array}
$$
\end{lemma}

{\bf Proof:} Suppose $g=\left (\begin{array}{cc}  a&b\\c&d\\ \end{array}\right  ) $ and consider the path $\gamma_0$. The first edge in the path is of type $A$ and is oriented forward, the second is type $C$ oriented backwards, and the final edge is type $A$ oriented backwards. The matrix $g$ is used to determine the contribution of the first two edges, while $\left (\begin{array}{cc}  a&a+b\\c&c+d\\ \end{array}\right  )$ is used for the third edge. Thus to determine the contribution of the third edge we must consider $-d^\prime/c^\prime=-(c+d)/c=-1-d/c$. Suppose first that $-d/c=\pm \infty$ and hence $-d^\prime/c^\prime=\pm \infty$. Using Table~\ref{table:m1andm2}, we see that the three edges contribute $(0,0), (0,-2\beta)$ and $(0,0)$ respectively, for a sum of $(0,-2\beta)$.
Next, suppose $-\infty<-d/c<0$ and thus  $-\infty<-d^\prime/c^\prime<0$. Now the edges contribute $(\beta, \beta), (0,-2\beta)$ and $(-\beta, -\beta)$, giving the same total as before. As our next case, suppose $0<-d/c<1$ and hence $-\infty<-d^\prime/c^\prime<0$.  The edges now contribute $(-\beta,-\beta), (2\beta,0)$ and $(-\beta, -\beta)$, again giving a sum of $(0,-2\beta)$. Finally, if $1<-d/c<\infty$, then $0<-d'/c'<\infty$ and we obtain $(-\beta,-\beta)+(0,-2\beta)+(\beta, \beta)=(0,-2\beta)$.

The computations for the other regions are similar and are left to the reader.
\hfill $\Box$  

{\bf Proof of Theorem~\ref{thm:xyz}:} Let $\gamma$ be a path in $D_t$ with endpoints in $D_1$ and let $\gamma'$ be obtained from $\gamma$ by pushing each ``diagonal'' edge across a region of type $R_0$ through $R_3$ (as in Lemma~\ref{lemma:pathmovest}) to eliminate all edges of type $C$ and $D$. Thus $M(\gamma')=k(\alpha, \beta)$ for some integer $k$. Notice that any path in $D_1$ must pass through fractions whose denominators alternate in parity. Thus if $\gamma'$ starts at $\frac{1}{0}$ and ends at $\frac{p}{q}$ it must have a number of edges equivalent to $q \mod 2$. Since the first edge contributes zero to $k$ and every other edge contributes $\pm 1$, we have that $k$ and $q$ have opposite parity.

To go back to $\gamma$ from $\gamma'$ suppose that we must move across $n_0^+$ regions of type $R_0$ (or $R_2$) in the positive sense and $n_0^-$ in the negative sense. Similarly, let $n_1^+$ and $n_1^-$ be the number of regions of type $R_1$ (or $R_3$) that we must push across in the positive or negative sense respectively. Then
$$\begin{array}{rl}
M(\gamma)&=k(\alpha, \beta)+(n_0^+-n_0^-)(0, -2\beta)+(n_1^+-n_1^-)(-\alpha+\beta, \alpha-\beta)\\
&=((k-n_1^++n_1^-)\alpha+(n_1^+-n_1^-)\beta, (n_1^+-n_1^-)\alpha+(k-n_1^++n_1^--2n_0^++2n_0^-)\beta).
\end{array}
$$
Thus
$$\begin{array}{rl}
x&=k-n_1^++n_1^-\\
y&=n_1^+-n_1^-\\
z&=k-n_1^++n_1^--2n_0^++2n_0^-
\end{array}
$$
and $x \equiv z \mod 2$. Furthermore, $x+y=k$. Thus if the path begins at $\frac{1}{0}$ and ends at $\frac{p}{q}$ we see that $x+y \equiv 1+q \mod 2$. 
\hfill $\Box$

The following lemma, which is analogous to Lemma~\ref{lemma:pathmovest}, allows us to handle the  case where $C$-type edges are replaced with $A$-type edges in $D_1$.
The proof is similar to the proof of Lemma~\ref{lemma:pathmovest} and is left to the reader.

\begin{figure}
\psfrag{a}{$\frac{a+b}{c+d}$}
\psfrag{b}{$\frac{a+2 b}{c+2 d}$}
\psfrag{c}{$\frac{b}{d}$}
\psfrag{d}{$\frac{a}{c}$}

\psfrag{e}{A}   
\psfrag{f}{B}   
\psfrag{g}{C}   
\psfrag{h}{D} 
\psfrag{e0}{A}   
\psfrag{f0}{B}   
\psfrag{g0}{C}   
\psfrag{h0}{D}
\psfrag{r}{D$_1$}

\psfrag{r0}{$S_0$}
\psfrag{r1}{$S_1$}
  
\begin{center}
    \leavevmode
    \scalebox{.85}{\includegraphics{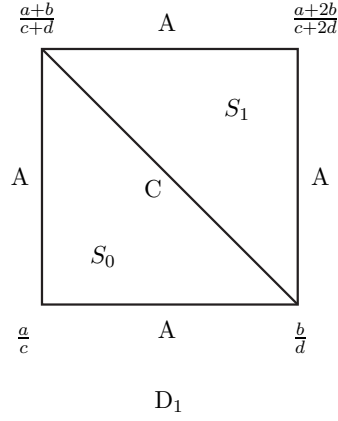}}
    \end{center}
\caption{Lemma~\ref{lemma:pathmoves1} describes the result of travelling around each region.}
\label{fig:gQ1}
\end{figure}

\begin{lemma}
\label{lemma:pathmoves1}
 Let  $g=\left (\begin{array}{cc}  a&b\\c&d\\ \end{array}\right  )$  be any element of $G$ where $b$ is even and label the regions of $g(Q)$ as shown in Figure~\ref{fig:gQ1}. If $\gamma_i=\partial S_i$, oriented counterclockwise, then
$$
\begin{array}{rl}
M(\gamma_0)&=( -2\beta+2n, -2n) \mbox{ and}\\
M(\gamma_1)&=(-2n,  -2\beta+2n)
\end{array}
$$
where $0 \le n \le \beta$.

\end{lemma}
 
Using this lemma we may determine the form of $M(\gamma)$ for any path in $D_1$.

\begin{theorem} 
\label{thm:xys}Let $\gamma$ be any path in $D_1$. Then $M(\gamma)$ has the form
$$M(\gamma)=((x+y s)\beta, (x-y s)\beta)$$
where $y$ is the number of $C$-type edges in $\gamma$ and $s$ is a rational parameter with $-1\le s \le 1$. If additionally, $\gamma$ begins at $\frac{1}{0}$ and ends at $\frac{p}{q}$, then $x+y \equiv 1+q \mod 2$.
\end{theorem}

{\bf Proof:}  If $\gamma$ has no $C$-type edges the result is the same as Lemma~\ref{lemma:AandBonly}. If $C$-type edges are present, push each such edge across a region of type $S_0$ (as in Lemma~\ref{lemma:pathmoves1}) to obtain a path with only $A$ and $B$-type edges. Applying Lemma~\ref{lemma:pathmoves1}, we now have that
$$M(\gamma)=k(\beta, \beta)+\sum_{i=1}^P(-2\beta+2n_i, -2n_i)-\sum_{i=P+1}^{P+N}(-2\beta+2n_i, -2n_i)$$
for some integer $k$, nonnegative integers $P$ and $N$, and integers $n_i$ where $0\le n_i \le \beta$. Let 
$$X=2\sum_{i=1}^Pn_i-2\sum_{i=P+1}^{P+N}n_i.$$
 If we let
$$s=\frac{X-(P-N)\beta}{(P+N)\beta}$$
then it is not hard to show that $-1 \le s \le 1$. Substituting for $X$ in $M(\gamma)$ we obtain the desired result with $y=P+N$ and $x=k-P+N$.
\hfill $\Box$

Theorems~\ref{thm:xyz} and \ref{thm:xys} place restrictions on the form of $M(\gamma)$ where $\gamma$ is a path in either $D_t$ or $D_1$. The following theorem shows that these are in fact the only restrictions that apply and moreover that minimal paths may be used to realize any desired value of $M(\gamma)$.

\begin{theorem}
There exist minimal paths in either $D_t$ or $D_1$ which begin at $\frac{1}{0}$ and end at some $\frac{p}{q}$ that realize all possible values of $M$ subject only to the constraints of Theorems~\ref{thm:xyz} and \ref{thm:xys}.
\end{theorem}

{\bf Proof:} Focussing on Theorem~\ref{thm:xyz}, let $x, y$, and $z$ be any three integers such that $x \equiv z \mod 2$. Suppose also that $x+y \equiv 1 \mod 2$. We seek a fraction $\frac{p}{q}$, with $q$ even, and a minimal path $\gamma$ from $\frac{1}{0}$ to $\frac{p}{q}$ such that $M(\gamma)=(x\alpha+y\beta, y\alpha+z\beta)$. Referring to the proof of Theorem~\ref{thm:xyz}, we will show how to build a minimal path with compete control over $k, n_0^+, n_0^-, n_1^+$, and $n_1^-$.

\begin{figure}
\psfrag{a}{$1$}
\psfrag{b}{$0$}
\psfrag{c}{$\frac{1}{0}$}
\psfrag{d}{$-1$}
\psfrag{e}{$\begin{array}{l}\Delta n_1^+=1\\ \Delta k=2 \end{array}$}
\psfrag{g}{$\begin{array}{l}\Delta n_1^-=1\\ \Delta k=0 \end{array}$}
\psfrag{i}{$\begin{array}{l}\Delta n_0^+=1\\ \Delta k=4 \end{array}$}
\psfrag{l}{$\begin{array}{l}\Delta n_0^-=1\\ \Delta k=0 \end{array}$}

\psfrag{n}{$\Delta k=1, 3, 5, \dots$}
\psfrag{p}{$\Delta k=-1,-3,-5, \dots$}

    \begin{center}
    \leavevmode
    \scalebox{.5}{\includegraphics{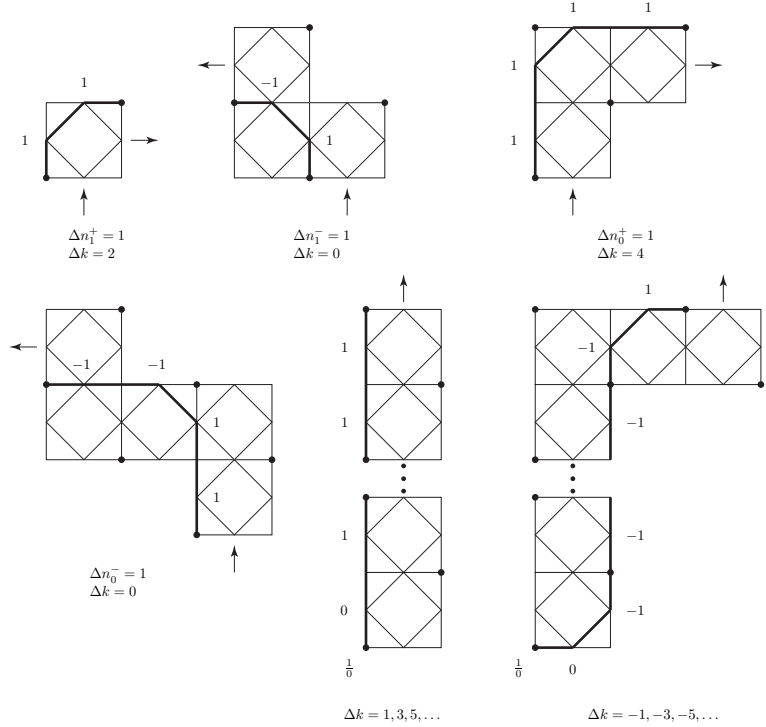}}
    \end{center}
\caption{Building blocks for a minimal path $\gamma$ with desired $M(\gamma)$.}
\label{fig:pathblocks}
\end{figure}

Figure~\ref{fig:pathblocks} shows six blocks of quadrilaterals, each containing a minimal path. Heavy dots are placed at vertices with even denominators. Starting with either of the last two blocks in the second row,  we may then paste on any number of the first four blocks (in any order), always gluing blocks together as indicated by the arrows.  No matter how  the blocks are joined  together, the path will remain minimal.  The first four blocks involve $C$ and $D$ type edges and can be used to create any desired values for  $n_0^+, n_0^-, n_1^+$, and $n_1^-$. Once these parameters are fixed, a sufficiently long starting block of one of the two types may be prepended to create any desired value of $k$.

It is not hard to adapt this argument to the case where $x+y \equiv 0 \mod 2$ or to the case of Theorem~\ref{thm:xys}. Moreover, more efficient sets of building blocks than these can easily be designed.
\hfill $\Box$

\section{An Example}
\label{section:example}

\begin{figure}[t]
\psfrag{a}{\large $\frac{1}{k}$}
    \begin{center}
    \leavevmode
        \scalebox{.85}{\includegraphics{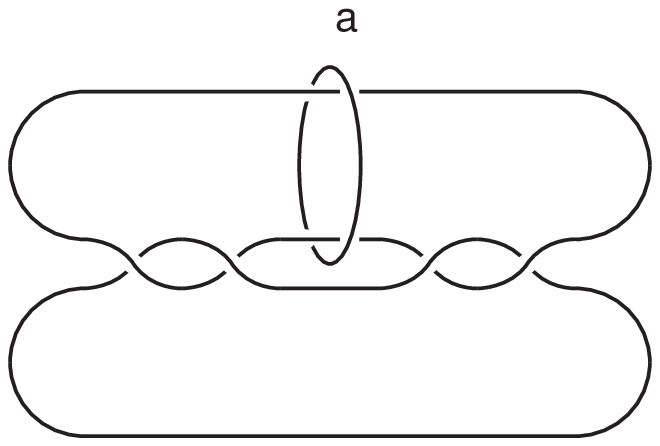}}
    \end{center}
\caption{The links ${\cal L}_k$}
\label{figure:the links}
\end{figure}

In this section we will apply our improved algorithm to compute the boundary slopes of the links
 ${\cal L}_\frac{4 k-1}{8 k}$. Since $\frac{4 k-1}{8 k}$ has the continued fraction expansion $[0,2,-2k,-2]$, these links can be pictured as shown in Figure~\ref{figure:the links}. Here we have replaced the $-2k$ right-handed crossings with $\frac{1}{k}$-surgery on an unknot surrounding the parallel strands. Viewed this way we see that ${\cal L}_\frac{4 k-1}{8 k}$, which we will henceforth simply denote  as ${\cal L}_k$,  is $\frac{1}{k}$ surgery on one component of the Borromean rings.

\begin{figure}[b]
\psfrag{a}{$\frac{1}{0}$}
\psfrag{b}{$S_0$}
\psfrag{c}{$\frac{1}{1}$}
\psfrag{d}{$\frac{2k}{4k+1}$}
\psfrag{e}{$T_{2k}$}
\psfrag{f}{$\frac{4k-1}{8k}$}
\psfrag{h}{$T_0$}
\psfrag{i}{$R_0$}
\psfrag{j}{$\frac{1}{2}$}
\psfrag{k}{$R_{2k+1}$}
\psfrag{l}{$S_{2k}$}
\psfrag{m}{$\frac{0}{1}$}
\psfrag{n}{$R_1$}
\psfrag{o}{$R_{2k}$}
\psfrag{p}{$\frac{2k-1}{4k-1}$}
\psfrag{q}{$R_i$}
\psfrag{r}{$R_{i+1}$}
\psfrag{s}{$\frac{i-1}{2 i-1}$}
\psfrag{t}{$\frac{i}{2i+1}$}
\psfrag{u}{$S_i$}
\psfrag{v}{$T_i$}
\psfrag{w}{$\frac{2i-1}{4i}$}

    \begin{center}
   \leavevmode
      \scalebox{.7}{\includegraphics{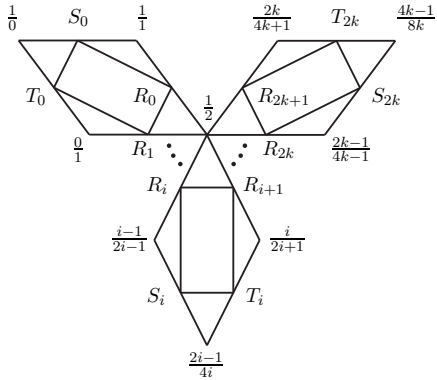}}
    \end{center}
\caption{The minimal chain of quadralaterals from $\frac{1}{0}$ to $\frac{4k-1}{8k}$.}
\label{figure:quadchain}
\end{figure}

The minimal chain of quadralaterals in the diagram $D_t$ connecting $\frac{1}{0}$ to $\frac{4k-1}{8k}$ is abstractly depicted in Figure~\ref{figure:quadchain}. It is not hard to show by induction that for $k>1$ there are precisely six minimal paths, $\{\gamma_1, \dots, \gamma_6 \}$, in $D_t$ from $\frac{1}{0}$ to $\frac{4k-1}{8k}$. These are listed in Table~\ref{table:minimalpaths} by listing the consecutive vertices in each path. The vertices $R_i, S_i$, and $T_i$ each lie between two vertices of $D_1$ (or $D_0$) as indicated in the figure. If $k=1$, path $\gamma_4$ is no longer minimal and should be deleted from the list.

Associated to each of these paths, in the case where $t>1$, is a weighted branched surface with $\alpha>\beta$. To find the boundary slopes of surfaces carried by such a branched surface we must first compute the intersection numbers $(M_1, M_2)$ and then adjust for our choice of basis. We will derive  $(M_1, M_2)$ from each path by using the Lemmas given in the last section.

\begin{table}
\renewcommand{\arraystretch}{1.5}
\caption{Minimal $D_t$ paths from $\frac{1}{0}$ to $\frac{4 k-1}{8k}$.}
\label{table:minimalpaths}
\begin{center}
\begin{tabular}{rl}
$\gamma_1$: & $\frac{1}{0}, S_0, R_0, \frac{1}{2}, R_{2 k+1}, T_{2k}, \frac{4k-1}{8k}$\\ 
$\gamma_2$: & $\frac{1}{0}, S_0, R_0, \frac{1}{2}, R_{2 k},  S_{2k}, \frac{4k-1}{8k}$\\ 
$\gamma_3$: &$\frac{1}{0}, T_0, R_1, \frac{1}{2}, R_{2 k+1},  T_{2k}, \frac{4k-1}{8k}$\\ 
$\gamma_4$: & $\frac{1}{0}, T_0, R_1, \frac{1}{2}, R_{2 k},  S_{2k}, \frac{4k-1}{8k}$\\ $\gamma_5$: & $\frac{1}{0}, T_0, R_1, R_2, \dots, R_{2k}, S_{2k}, \frac{4k-1}{8k}$\\ 
$\gamma_6$:  & $\frac{1}{0}, T_0, \frac{0}{1}, S_1, T_1, \frac{1}{3}, S_2, T_2, \dots, \frac{2k-1}{4k-1}, S_{2k},  \frac{4k-1}{8k}$

\end{tabular}
\end{center}
\end{table}

Let $\sigma_1$ be the path in $D_1$ given by $\sigma_1=\{\frac{1}{0},\frac{1}{1},\frac{1}{2},\frac{2k}{4k+1},\frac{4k-1}{8k}\}$. It follows easily from Lemma~\ref{lemma:AandBonly} that $M(\sigma_1)=(3\alpha, 3\beta)$. We may now deform $\sigma_1$ to $\gamma_1$ by moving the path over two regions, each of type $R_1$ (or $R_3$) as shown in Figure~\ref{fig:gQ}. This results in adding $2(-\alpha+\beta, \alpha-\beta)$. We therefore obtain $(2 \alpha+\beta, \alpha+2 \beta)$.
To pass to either $\gamma_2$ or $\gamma_3$, we deform $\gamma_1$ across two regions of type $R_0$ and one of type $R_4$. In either case then, we must add 
$(-2\beta, -2\alpha)$ and obtain a total of $(\alpha, \beta)$.
To move from $\gamma_2$ (or $\gamma_3$) to $\gamma_4$, we again move across two regions of type $R_0$ and one of type $R_4$ and again add  $(-2\beta, -2\alpha)$. This gives $(\alpha-2\beta, -2\alpha+\beta)$ for $M(\gamma_4)$.
We may now obtain $\gamma_5$ from $\gamma_4$ by moving across $2k-1$ regions of type $R_0$. Thus the value of $M(\gamma_5)$ is $(\alpha-2\beta, -2\alpha+\beta)+(2k-1)(0, -2\beta)=(\alpha-2\beta, -2\alpha+(3-4k)\beta)$.
Now let $\sigma_2=\{\frac{1}{0},\frac{0}{1},\frac{1}{2},\frac{2k-1}{4k-1},\frac{4k-1}{8k}\}$. We may either think of moving $\gamma_4$ to $\sigma_2$, or start over again with Lemma~\ref{lemma:pathmovest}. Either way we obtain $(-\alpha, -\beta)$ for $M(\sigma_2)$.
Finally, to obtain $\gamma_6$ from $\sigma_2$, we must move across $2k-1$ regions, each made up of one region of type $R_0$, two of type $R_1$ and one of type $R_4$. This results in adding $(2k-1)(-2\alpha, 0)$ to  $M(\sigma_2)$ giving a total of $((1-4k)\alpha, -\beta)$.

We have now accounted for all the possibilities when $\alpha>\beta$. We may determine all cases where $\beta=0$, and thus $t=\infty$, from these by substitution, but letting $\alpha$ equal $\beta$  does not give all possibilities for $t=1$. Because both $\gamma_5$ and $\gamma_6$ involve $C$-type edges, so do their limiting minimal edge paths in $D_1$.  In this example, both $\gamma_5$ and $\gamma_6$  limit to the same minimal edge path $\gamma_5^1$ in $D_1$. To recover all possible boundary slopes when $\alpha=\beta$ we must now consider the branched surface corresponding to this path where $C$-type surfaces are replaced with $A$-type surfaces.

Starting from $\sigma_2$ we may move to $\gamma_5^1$ by moving across $k$ regions of type $S_0$ and $k-1$ regions of type $S_1$. Thus we must add $\sum_{i=1}^k(-2\beta+2n_{2i-1}, -2n_{2i-1})$ and $\sum_{i=1}^{k-1}(-2n_{2i}, -2 \beta+2n_{2i})$ to $M(\sigma_2)=(-\alpha, -\beta)$. This gives $(-(1+2k)\beta+X, (1-2k)\beta-X)$, where $X=2\sum_{i=1}^{2k-1}(-1)^{i+1}n_i$. Here $0\le n_i\le \beta$.

The final 4-tuples $(M_1, \alpha, M_2, \beta)$ for ${\cal L}_{\frac{4k-1}{8k}}$ in all cases where $\alpha \ge \beta \ge 0$ are given in  Table~\ref{table:boundaryslopes}. All of these data  are still with respect to the basis $\{\mu_i, \lambda_i\}$. 

From Table~\ref{table:boundaryslopes} we may easily derive the corresponding 4-tuples when $0 \le \alpha \le \beta$. Namely, $(M_1(\alpha, \beta),\alpha, M_2(\alpha, \beta), \beta)$ is changed to $(M_2(\beta, \alpha), \alpha, M_1(\beta, \alpha), \beta)$. Notice that if  a minimal path in $D_t$ does not involve any $C$-type edges, then the contributions to $M_1$ and $M_2$ given in Table~\ref{table:m1andm2} are equal when $\alpha=\beta$. Thus the same 4-tuple will result as we approach $\frac{\alpha}{\beta}=1$ from either above or below.
Using Equation~\ref{eqn:linkingnumber} it is not difficult to show that $\mbox{lk}(\lambda_1, K_1)=\mbox{lk}(\lambda_2, K_2)=-1$ for all $k$.  The final boundary slopes $M_1/\alpha$ and $M_2/\beta$  for all rational values of $t$, and with respect to the preferred basis are now given in Table~\ref{table:finalslopes}. Furthermore, the slopes in  this table are described by the rational parameters $t=\alpha/\beta$ in the case where $\alpha\ne \beta$, and  
 $$s=\frac{X-\beta}{(2k-1)\beta}$$ 
  in the case where $\alpha=\beta$. Note that a pair of boundary slopes of the form $(0,\phi)$ means that the corresponding surface has no boundary components on $K_2$ (and a slope of zero on $K_1$).

\begin{table}
\caption{Boundary slope data for ${\cal L}_{\frac{4k-1}{8k}}$, $k>0$,  with respect to the basis $\{\mu_1, \lambda_1\}$ and for $\alpha \ge \beta$.
The path $\gamma_4$ is not counted if $k=1$. Here $X=n_1-n_2+n_3- \dots +n_{2k-1}$.}
\begin{center}
\begin{tabular}{lcc}
path&algebraic intersection with & restrictions\\
&$(\lambda_1, \mu_1, \lambda_2, \mu_2)$&\\
\hline
$\gamma_1$ &$(\alpha+2 \beta,\alpha,2 \alpha+\beta, \beta)$&$\alpha>\beta \ge 0$\\
$\gamma_2$ &$(\alpha,\alpha,\beta, \beta)$&$\alpha>\beta\ge 0$\\
$\gamma_3$ &$(\alpha,\alpha,\beta, \beta)$&$\alpha>\beta\ge 0$\\
$\gamma_4$ &$(\alpha-2 \beta,\alpha,-2\alpha+\beta, \beta)$&$\alpha>\beta\ge 0$\\
$\gamma_5$ &$(\alpha-2 \beta,\alpha,-2 \alpha+(3-4k)\beta, \beta)$&$\alpha>\beta\ge 0$\\
$\gamma_6$ &$((1-4k)\alpha,\alpha, -\beta, \beta)$&$\alpha>\beta\ge 0$\\
$\gamma_5^1$ &$(-(1+2k)\beta+2X,\beta, (1-2k)\beta-2X, \beta)$& $0\le n_i\le \beta=\alpha$\\

\end{tabular}
\end{center}
\label{table:boundaryslopes}
\end{table}

\begin{table}
\caption{Boundary slope pairs for  ${\cal L}_{\frac{4k-1}{8k}}$, $k>0$,  with respect to the preferred basis $\{\mu_1, \lambda_1^0\}$. Both $t$ and $s$ are rational parameters.}
\begin{center}
\begin{tabular}{|l|c|}
\hline
$\partial$-slopes&restrictions\\
\hline
\hline
$(0,0)$&\\
$(0, \phi),(\phi, 0)$&\\
$(-4k, \phi),(\phi, -4k)$&\\
$(-4k, -2),(-2, -4k)$&\\
$(2 t^{-1}, 2t)$&$0\le t \le \infty$\\
$(-2 t^{-1}, -2t)$&$0\le t \le \infty, k>1$\\
$(-2 t^{-1}+2-4k, -2t)$&$0\le t \le 1$\\
$(-2 t^{-1},2-4k- 2t)$&$1\le t \le \infty$\\
$(-1-2k+(2k-1)s, -1-2k-(2k-1)s)$&$-1 \le s\le 1$\\

\hline
\end{tabular}
\end{center}
\label{table:finalslopes}
\end{table}

\section{Boundary Slopes of 2-bridge Links up to 10 Crossings}
\label{section:table}

We have written a computer program to implement the algorithm described in Section~\ref{section:improvedalgorithm}. One begins by finding all minimal paths in $D_0$ and $D_1$ from $\frac{1}{0}$ to $\frac{p}{q}$. These determine the minimal paths in $D_t$. Each minimal path $\gamma$ in $D_t$ is then deformed into a path $\gamma'$ of only $A$ and $B$-type edges. Lemma~\ref{lemma:AandBonly} is used to compute $M(\gamma')$ and $M(\gamma)$ is derived from this using Lemma~\ref{lemma:pathmovest}.  In order to minimize the chance of producing errors, each of us coded a program independent of the other with debugging proceeding until our results agreed.

The second Tait conjecture states that any reduced alternating diagram of a link is minimal. Thus is it easy to determine the crossing number of ${\cal L}_{p/q}$. If we express $p/q$ as the continued fraction $[0, a_2, \dots, a_n]$ where each $a_i$ is positive, then ${\cal L}_{p/q}$ has a reduced alternating 4-plat diagram with $a_2+a_3+\dots a_n$ crossings. Thus the crossing number of ${\cal L}_{p/q}$ is $a_2+a_3+\dots a_n$. It is now a simple matter to determine all 2-bridge links with 10 or fewer crossings by finding all such continued fractions. Furthermore, from the classification of 2-bridge links, we know precisely when two fractions represent the same link.

The boundary slopes of all 2-bridge links having ten or less crossings are presented in Tables~\ref{data2to8}-\ref{data10c}. All slopes are  with respect to a preferred longitude and meridian basis. For those links with 9 or less crossings, the index of the link in Rolfsen's table \cite{Ro:1990} is also given. Two types of data are listed. First, all boundary slopes in the case where $1<t<\infty$  are given. From these the slopes with $0<t<1$ can then be determined as described earlier. Furthermore, one may then derive from these all boundary slopes in the cases $t=0$ and $t=\infty$ and some of the boundary slopes in the case $t=1$. Finally, all additional slopes in the case $t=1$ are given in terms of the rational parameter $s$ with $-1 \le s \le 1$.

\bibliographystyle{hplain}
%\bibliography{hs}

\begin{table}[htdp]
\caption{Boundary slope data for 2-bridge links to 8 crossings. The parameters $t$ and $s$ are rational with $1<t<\infty$ and $-1\le s \le 1$.}
\begin{center}

$\begin{array}{|c|r|llll|}
\hline
\mbox{link}&p/q&\mbox{boundary slopes}&&&\\
\hline
\hline
2^2_{1}&1/2&(-t^{-1},-t)&(t^{-1},t)&&\\
\hline
4^2_{1}&1/4&(-2,-2)&(-2t^{-1},-2 t)&(2t^{-1},2 t)&\\
\hline
5^2_{1}&3/8&(-4,-2)&(0,0)&(-2t^{-1},-2 - 2 t)&(2t^{-1},2 t)\\
&&(-3 + s,-3 - s)&&&\\
\hline
6^2_{1}&1/6&(-3,-3)&(-3t^{-1},-3 t)&(3t^{-1},3 t)&\\
6^2_{2}&3/10&(-2-t^{-1},-2 - t)&(2-t^{-1},-t)&(-2+t^{-1},t)&(2+t^{-1},2 + t)\\
&&(-3t^{-1},-3 t)&(3t^{-1},3 t)&(s,-s)&\\
6^2_{3}&5/12&(-6,-2)&(0,0)&(-2t^{-1},-4 - 2 t)&(-2t^{-1},-2 t)\\
&&(2t^{-1},2 t)&(-4 + 2 s,-4 - 2 s)&&\\
\hline
7^2_{1}&3/14&(-5,-3)&(-2-t^{-1},-4 - t)&(-2+t^{-1},-2 + t)&(-3t^{-1},-2 - 3 t)\\
&&(-t^{-1},-t)&(t^{-1},t)&(3t^{-1},3 t)&(-4 + s,-4 - s)\\
7^2_{3}&7/16&(-8,-2)&(0,0)&(-2t^{-1},-6 - 2 t)&(-2t^{-1},-2 t)\\
&&(2t^{-1},2 t)&(-5 + 3 s,-5 - 3 s)&&\\
7^2_{2}&5/18&(-2-t^{-1},-2 - t)&(4-t^{-1},-t)&(-2+t^{-1},-2 + t)&(-2+t^{-1},2 + t)\\
&&(4+t^{-1},2 + t)&(-3t^{-1},-3 t)&(-t^{-1},-t)&(t^{-1},t)\\
&&(3t^{-1},2 + 3 t)&(4 + s,4 - s)&(1 + 2 s,1 - 2 s)&\\
\hline
8^2_{1}&1/8&(-4,-4)&(-4t^{-1},-4 t)&(4t^{-1},4 t)&\\
8^2_{2}&3/16&(0,-2)&(2-2t^{-1},-2 t)&(-3-t^{-1},-3 - t)&(-3+t^{-1},-1 + t)\\
&&(2+2t^{-1},2 + 2 t)&(-4t^{-1},-4 t)&(4t^{-1},4 t)&(-1 + s,-1 - s)\\
8^2_{6}&9/20&(-10,-2)&(0,0)&(-2t^{-1},-8 - 2 t)&(-2t^{-1},-2 t)\\
&&(2t^{-1},2 t)&(-6 + 4 s,-6 - 4 s)&&\\
8^2_{3}&5/22&(-7,-3)&(-2-t^{-1},-6 - t)&(-2-t^{-1},-2 - t)&(-2+t^{-1},-2 + t)\\
&&(-3t^{-1},-4 - 3 t)&(-3t^{-1},-3 t)&(-t^{-1},-t)&(t^{-1},t)\\
&&(3t^{-1},3 t)&(-5 + 2 s,-5 - 2 s)&&\\
8^2_{4}&7/24&(-4,-4)&(-4,0)&(2,0)&(-2-2t^{-1},-2 - 2 t)\\
&&(2-2t^{-1},-2 - 2 t)&(-2+2t^{-1},2 t)&(2+2t^{-1},2 + 2 t)&(-4t^{-1},-4 t)\\
&&(4t^{-1},4 t)&(1 + s,1 - s)&(-2 + 2 s,-2 - 2 s)&\\
8^2_{5}&7/26&(-2-t^{-1},-2 - t)&(6-t^{-1},-t)&(-2+t^{-1},-2 + t)&(-2+t^{-1},4 + t)\\
&&(6+t^{-1},2 + t)&(-3t^{-1},-3 t)&(-t^{-1},-t)&(t^{-1},t)\\
&&(3t^{-1},3 t)&(3t^{-1},4 + 3 t)&(5 + 2 s,5 - 2 s)&(2 + 3 s,2 - 3 s)\\
8^2_{7}&11/30&(-7,-3)&(-4-t^{-1},-4 - t)&(-4+t^{-1},-2 + t)&(-3t^{-1},-4 - 3 t)\\
&&(-t^{-1},-2 - t)&(-t^{-1},-t)&(t^{-1},t)&(3t^{-1},3 t)\\
&&(-2 + s,-2 - s)&(-5 + 2 s,-5 - 2 s)&&\\
8^2_{8}&13/34&(-4-t^{-1},-2 - t)&(4-t^{-1},-2 - t)&(4-t^{-1},2 - t)&(-4+t^{-1},-2 + t)\\
&&(-4+t^{-1},2 + t)&(4+t^{-1},2 + t)&(-3t^{-1},-2 - 3 t)&(-t^{-1},-2 - t)\\
&&(-t^{-1},-t)&(t^{-1},t)&(t^{-1},2 + t)&(3t^{-1},2 + 3 t)\\
&&(-4 + s,-4 - s)&(-2 + s,-2 - s)&(3 s,-3 s)&(2 + s,2 - s)\\
&&(4 + s,4 - s)&&&\\
\hline \end{array} $
\end{center}
\label{data2to8}
\end{table}

\begin{table}[htdp]
\caption{Boundary slope data for 2-bridge links with 9 crossings. The parameters $t$ and $s$ are rational with $1<t<\infty$ and $-1 \le s \le 1$.}
\begin{center}

$\begin{array}{|c|r|llll|}
\hline
\mbox{link}&p/q&\mbox{boundary slopes}&&&\\
\hline
\hline
9^2_{1}&3/20&(-6,-4)&(-3-t^{-1},-5 - t)&(-3+t^{-1},-3 + t)&(-4t^{-1},-2 - 4 t)\\
&&(-2t^{-1},-2 t)&(2t^{-1},2 t)&(4t^{-1},4 t)&(-5 + s,-5 - s)\\
9^2_{4}&5/24&(-6,-4)&(-4,-6)&(0,0)&(-2-2t^{-1},-4 - 2 t)\\
&&(-2+2t^{-1},-2 + 2 t)&(-4t^{-1},-2 - 4 t)&(4t^{-1},4 t)&(-5 + s,-5 - s)\\
9^2_{10}&11/24&(-12,-2)&(0,0)&(-2t^{-1},-10 - 2 t)&(-2t^{-1},-2 t)\\
&&(2t^{-1},2 t)&(-7 + 5 s,-7 - 5 s)&&\\
9^2_{2}&5/28&(2,-2)&(4-2t^{-1},-2 t)&(-3-t^{-1},-3 - t)&(-3+t^{-1},-3 + t)\\
&&(-3+t^{-1},1 + t)&(4+2t^{-1},2 + 2 t)&(-4t^{-1},-4 t)&(-2t^{-1},-2 t)\\
&&(2t^{-1},2 t)&(4t^{-1},2 + 4 t)&(2 s,-2 s)&(5 + s,5 - s)\\
9^2_{3}&7/30&(-9,-3)&(-2-t^{-1},-8 - t)&(-2-t^{-1},-2 - t)&(-2+t^{-1},-2 + t)\\
&&(-3t^{-1},-6 - 3 t)&(-3t^{-1},-3 t)&(-t^{-1},-t)&(t^{-1},t)\\
&&(3t^{-1},3 t)&(-6 + 3 s,-6 - 3 s)&&\\
9^2_{5}&7/32&(0,-4)&(0,0)&(-2-2t^{-1},-4 - 2 t)&(2-2t^{-1},-2 - 2 t)\\
&&(2-2t^{-1},2 - 2 t)&(-5-t^{-1},-3 - t)&(-5+t^{-1},-1 + t)&(-2+2t^{-1},-2 + 2 t)\\
&&(2+2t^{-1},2 + 2 t)&(-4t^{-1},-2 - 4 t)&(4t^{-1},4 t)&(-5 + s,-5 - s)\\
&&(-2 + 2 s,-2 - 2 s)&&&\\
9^2_{8}&9/34&(-2-t^{-1},-2 - t)&(8-t^{-1},-t)&(-2+t^{-1},-2 + t)&(-2+t^{-1},6 + t)\\
&&(8+t^{-1},2 + t)&(-3t^{-1},-3 t)&(-t^{-1},-t)&(t^{-1},t)\\
&&(3t^{-1},3 t)&(3t^{-1},6 + 3 t)&(6 + 3 s,6 - 3 s)&(3 + 4 s,3 - 4 s)\\
9^2_{6}&11/36&(-2,-2)&(-2,0)&(2,2)&(-2-2t^{-1},-2 - 2 t)\\
&&(2-2t^{-1},-2 t)&(5-t^{-1},1 - t)&(5+t^{-1},3 + t)&(-2+2t^{-1},2 + 2 t)\\
&&(2+2t^{-1},4 + 2 t)&(-4t^{-1},-4 t)&(-2t^{-1},-2 t)&(2t^{-1},2 t)\\
&&(4t^{-1},2 + 4 t)&(-1 + s,-1 - s)&(5 + s,5 - s)&(2 + 2 s,2 - 2 s)\\
9^2_{9}&11/40&(-4,-4)&(-4,2)&(0,0)&(4,0)\\
&&(-2-2t^{-1},-2 - 2 t)&(4-2t^{-1},-2 - 2 t)&(-2+2t^{-1},-2 + 2 t)&(-2+2t^{-1},2 + 2 t)\\
&&(4+2t^{-1},2 + 2 t)&(-4t^{-1},-4 t)&(4t^{-1},2 + 4 t)&(5 + s,5 - s)\\
&&(2 + 2 s,2 - 2 s)&(-1 + 3 s,-1 - 3 s)&&\\
9^2_{7}&13/44&(-6,-4)&(-6,0)&(-2,-2)&(-2,0)\\
&&(2,-2)&(2,0)&(2,2)&(-4-2t^{-1},-2 - 2 t)\\
&&(-2-2t^{-1},-4 - 2 t)&(2-2t^{-1},-4 - 2 t)&(2-2t^{-1},-2 t)&(-4+2t^{-1},2 t)\\
&&(-2+2t^{-1},2 t)&(2+2t^{-1},2 + 2 t)&(-4t^{-1},-2 - 4 t)&(-2t^{-1},-2 t)\\
&&(2t^{-1},2 t)&(4t^{-1},4 t)&(-5 + s,-5 - s)&(-1 + s,-1 - s)\\
&&(2 s,-2 s)&(1 + s,1 - s)&(-3 + 3 s,-3 - 3 s)&\\
9^2_{11}&17/46&(-9,-3)&(-6-t^{-1},-4 - t)&(-4-t^{-1},-6 - t)&(-4-t^{-1},-2 - t)\\
&&(-6+t^{-1},-2 + t)&(-4+t^{-1},-2 + t)&(-3t^{-1},-6 - 3 t)&(-3t^{-1},-2 - 3 t)\\
&&(-t^{-1},-4 - t)&(-t^{-1},-2 - t)&(-t^{-1},-t)&(t^{-1},t)\\
&&(3t^{-1},3 t)&(-4 + s,-4 - s)&(-2 + s,-2 - s)&(-3 + 2 s,-3 - 2 s)\\
&&(-6 + 3 s,-6 - 3 s)&&&\\
9^2_{12}&19/50&(-4-t^{-1},-2 - t)&(6-t^{-1},-2 - t)&(6-t^{-1},2 - t)&(-4+t^{-1},-2 + t)\\
&&(-4+t^{-1},4 + t)&(6+t^{-1},2 + t)&(-3t^{-1},-2 - 3 t)&(-t^{-1},-2 - t)\\
&&(-t^{-1},-t)&(t^{-1},t)&(t^{-1},4 + t)&(3t^{-1},3 t)\\
&&(3t^{-1},4 + 3 t)&(-4 + s,-4 - s)&(-2 + s,-2 - s)&(3 + 2 s,3 - 2 s)\\
&&(5 + 2 s,5 - 2 s)&(1 + 4 s,1 - 4 s)&&\\
\hline \end{array} $
\end{center}
\label{data9}
\end{table}

\begin{table}[htdp]
\caption{Boundary slope data for 2-bridge links with 10 crossings (part 1). The parameters $t$ and $s$ are rational with $1<t<\infty$ and $-1 \le s \le 1$. }
\begin{center}

$\begin{array}{|r|llll|}
\hline
p/q&\mbox{boundary slopes}&&&\\
\hline
\hline
1/10&(-5,-5)&(-5t^{-1},-5 t)&(5t^{-1},5 t)&\\
3/22&(-1,-3)&(2-3t^{-1},-3 t)&(-4-t^{-1},-4 - t)&(-4+t^{-1},-2 + t)\\
&(2+3t^{-1},2 + 3 t)&(-5t^{-1},-5 t)&(5t^{-1},5 t)&(-2 + s,-2 - s)\\
5/26&(-1,1)&(1,-1)&(-3-2t^{-1},-3 - 2 t)&(3-2t^{-1},1 - 2 t)\\
&(-3+2t^{-1},-1 + 2 t)&(3+2t^{-1},3 + 2 t)&(-5t^{-1},-5 t)&(5t^{-1},5 t)\\
&(s,-s)&&&\\
13/28&(-14,-2)&(0,0)&(-2t^{-1},-12 - 2 t)&(-2t^{-1},-2 t)\\
&(2t^{-1},2 t)&(-8 + 6 s,-8 - 6 s)&&\\
5/32&(-8,-4)&(-3-t^{-1},-7 - t)&(-3-t^{-1},-3 - t)&(-3+t^{-1},-3 + t)\\
&(-4t^{-1},-4 - 4 t)&(-4t^{-1},-4 t)&(-2t^{-1},-2 t)&(2t^{-1},2 t)\\
&(4t^{-1},4 t)&(-6 + 2 s,-6 - 2 s)&&\\
7/38&(-5,-5)&(-5,-1)&(-2-3t^{-1},-2 - 3 t)&(2-3t^{-1},-2 - 3 t)\\
&(-3-2t^{-1},-3 - 2 t)&(2-t^{-1},-t)&(2+t^{-1},t)&(-3+2t^{-1},-1 + 2 t)\\
&(-2+3t^{-1},3 t)&(2+3t^{-1},2 + 3 t)&(-5t^{-1},-5 t)&(-t^{-1},-4 - t)\\
&(t^{-1},-2 + t)&(5t^{-1},5 t)&(s,-s)&(2 + s,2 - s)\\
&(-3 + 2 s,-3 - 2 s)&&&\\
9/38&(-11,-3)&(-2-t^{-1},-10 - t)&(-2-t^{-1},-2 - t)&(-2+t^{-1},-2 + t)\\
&(-3t^{-1},-8 - 3 t)&(-3t^{-1},-3 t)&(-t^{-1},-t)&(t^{-1},t)\\
&(3t^{-1},3 t)&(-7 + 4 s,-7 - 4 s)&&\\
7/40&(4,-2)&(6-2t^{-1},-2 t)&(-3-t^{-1},-3 - t)&(-3+t^{-1},-3 + t)\\
&(-3+t^{-1},3 + t)&(6+2t^{-1},2 + 2 t)&(-4t^{-1},-4 t)&(-2t^{-1},-2 t)\\
&(2t^{-1},2 t)&(4t^{-1},4 t)&(4t^{-1},4 + 4 t)&(6 + 2 s,6 - 2 s)\\
&(1 + 3 s,1 - 3 s)&&&\\
9/40&(-8,-4)&(-4,-8)&(-4,-4)&(0,0)\\
&(-2-2t^{-1},-6 - 2 t)&(-2-2t^{-1},-2 - 2 t)&(-2+2t^{-1},-2 + 2 t)&(-4t^{-1},-4 - 4 t)\\
&(-4t^{-1},-4 t)&(4t^{-1},4 t)&(-6 + 2 s,-6 - 2 s)&\\
11/42&(-2-t^{-1},-2 - t)&(10-t^{-1},-t)&(-2+t^{-1},-2 + t)&(-2+t^{-1},8 + t)\\
&(10+t^{-1},2 + t)&(-3t^{-1},-3 t)&(-t^{-1},-t)&(t^{-1},t)\\
&(3t^{-1},3 t)&(3t^{-1},8 + 3 t)&(7 + 4 s,7 - 4 s)&(4 + 5 s,4 - 5 s)\\
13/42&(3,1)&(-2-3t^{-1},-2 - 3 t)&(3-2t^{-1},-1 - 2 t)&(-4-t^{-1},-4 - t)\\
&(-4+t^{-1},t)&(3+2t^{-1},3 + 2 t)&(-2+3t^{-1},3 t)&(-5t^{-1},-5 t)\\
&(-t^{-1},-t)&(t^{-1},2 + t)&(5t^{-1},5 t)&(2 + s,2 - s)\\
&(-1 + 2 s,-1 - 2 s)&&&\\
11/48&(0,-6)&(0,0)&(-2-2t^{-1},-6 - 2 t)&(-2-2t^{-1},-2 - 2 t)\\
&(2-2t^{-1},-4 - 2 t)&(2-2t^{-1},2 - 2 t)&(-7-t^{-1},-3 - t)&(-7+t^{-1},-1 + t)\\
&(-2+2t^{-1},-2 + 2 t)&(2+2t^{-1},2 + 2 t)&(-4t^{-1},-4 - 4 t)&(-4t^{-1},-4 t)\\
&(4t^{-1},4 t)&(-6 + 2 s,-6 - 2 s)&(-3 + 3 s,-3 - 3 s)&\\
17/48&(-8,-4)&(-2,-4)&(0,0)&(-2-2t^{-1},-6 - 2 t)\\
&(-5-t^{-1},-5 - t)&(-5+t^{-1},-3 + t)&(-2+2t^{-1},-2 + 2 t)&(-4t^{-1},-4 - 4 t)\\
&(-2t^{-1},-2 - 2 t)&(2t^{-1},2 t)&(4t^{-1},4 t)&(-3 + s,-3 - s)\\
&(-6 + 2 s,-6 - 2 s)&&&\\
\hline \end{array} $

\end{center}
\label{data10a}
\end{table}

\begin{table}[htdp]
\caption{Boundary slope data for 2-bridge links with 10 crossings (part 2). The parameters $t$ and $s$ are rational with $1<t<\infty$ and $-1 \le s \le 1$.}
\begin{center}

$\begin{array}{|r|llll|}
\hline
p/q&\mbox{boundary slopes}&&&\\
\hline
\hline
11/52&(-8,-4)&(-6,-6)&(-2,-4)&(-2,-2)\\
&(0,-2)&(0,0)&(-4-2t^{-1},-4 - 2 t)&(-2-2t^{-1},-6 - 2 t)\\
&(-5-t^{-1},-5 - t)&(-5+t^{-1},-3 + t)&(-4+2t^{-1},-2 + 2 t)&(-2+2t^{-1},-2 + 2 t)\\
&(-4t^{-1},-4 - 4 t)&(-2t^{-1},-2 - 2 t)&(-2t^{-1},-2 t)&(2t^{-1},2 t)\\
&(4t^{-1},4 t)&(-3 + s,-3 - s)&(-1 + s,-1 - s)&(-6 + 2 s,-6 - 2 s)\\
15/56&(-4,-4)&(-4,4)&(0,0)&(6,0)\\
&(-2-2t^{-1},-2 - 2 t)&(6-2t^{-1},-2 - 2 t)&(-2+2t^{-1},-2 + 2 t)&(-2+2t^{-1},4 + 2 t)\\
&(6+2t^{-1},2 + 2 t)&(-4t^{-1},-4 t)&(4t^{-1},4 t)&(4t^{-1},4 + 4 t)\\
&(4 s,-4 s)&(6 + 2 s,6 - 2 s)&(3 + 3 s,3 - 3 s)&\\
17/56&(-2,-2)&(-2,0)&(2,0)&(2,2)\\
&(2,4)&(4,2)&(-2-2t^{-1},-2 - 2 t)&(2-2t^{-1},-2 t)\\
&(7-t^{-1},1 - t)&(7+t^{-1},3 + t)&(-2+2t^{-1},2 t)&(-2+2t^{-1},4 + 2 t)\\
&(2+2t^{-1},2 + 2 t)&(2+2t^{-1},6 + 2 t)&(-4t^{-1},-4 t)&(-2t^{-1},-2 t)\\
&(2t^{-1},2 t)&(4t^{-1},4 t)&(4t^{-1},4 + 4 t)&(-1 + s,-1 - s)\\
&(1 + s,1 - s)&(6 + 2 s,6 - 2 s)&(3 + 3 s,3 - 3 s)&\\
17/58&(-2-3t^{-1},-2 - 3 t)&(2-3t^{-1},-2 - 3 t)&(2-3t^{-1},-3 t)&(-4-t^{-1},-4 - t)\\
&(-4-t^{-1},-t)&(-2-t^{-1},-t)&(2-t^{-1},-t)&(4-t^{-1},-2 - t)\\
&(4-t^{-1},2 - t)&(-4+t^{-1},-2 + t)&(-4+t^{-1},2 + t)&(-2+t^{-1},t)\\
&(2+t^{-1},t)&(4+t^{-1},t)&(4+t^{-1},4 + t)&(-2+3t^{-1},3 t)\\
&(-2+3t^{-1},2 + 3 t)&(2+3t^{-1},2 + 3 t)&(-5t^{-1},-5 t)&(-t^{-1},-2 - t)\\
&(t^{-1},2 + t)&(5t^{-1},5 t)&(-2 + s,-2 - s)&(s,-s)\\
&(3 s,-3 s)&(2 + s,2 - s)&(-3 + 2 s,-3 - 2 s)&(3 + 2 s,3 - 2 s)\\
13/60&(-2,-4)&(-2,-2)&(0,0)&(0,2)\\
&(2,-4)&(2,0)&(-2-2t^{-1},-4 - 2 t)&(4-2t^{-1},-2 - 2 t)\\
&(4-2t^{-1},2 - 2 t)&(-5-t^{-1},-3 - t)&(-5+t^{-1},-3 + t)&(-5+t^{-1},1 + t)\\
&(-2+2t^{-1},2 t)&(4+2t^{-1},2 + 2 t)&(-4t^{-1},-2 - 4 t)&(-2t^{-1},-2 - 2 t)\\
&(-2t^{-1},-2 t)&(2t^{-1},2 t)&(4t^{-1},2 + 4 t)&(-5 + s,-5 - s)\\
&(-3 + s,-3 - s)&(1 + s,1 - s)&(5 + s,5 - s)&(-1 + 3 s,-1 - 3 s)\\
23/62&(-11,-3)&(-8-t^{-1},-4 - t)&(-4-t^{-1},-8 - t)&(-4-t^{-1},-2 - t)\\
&(-8+t^{-1},-2 + t)&(-4+t^{-1},-2 + t)&(-3t^{-1},-8 - 3 t)&(-3t^{-1},-2 - 3 t)\\
&(-t^{-1},-6 - t)&(-t^{-1},-2 - t)&(-t^{-1},-t)&(t^{-1},t)\\
&(3t^{-1},3 t)&(-4 + s,-4 - s)&(-2 + s,-2 - s)&(-4 + 3 s,-4 - 3 s)\\
&(-7 + 4 s,-7 - 4 s)&&&\\
19/64&(-8,-4)&(-8,0)&(-2,-2)&(-2,0)\\
&(2,-4)&(2,0)&(2,2)&(-6-2t^{-1},-2 - 2 t)\\
&(-2-2t^{-1},-6 - 2 t)&(-2-2t^{-1},-2 - 2 t)&(2-2t^{-1},-6 - 2 t)&(2-2t^{-1},-2 t)\\
&(-6+2t^{-1},2 t)&(-2+2t^{-1},2 t)&(2+2t^{-1},2 + 2 t)&(-4t^{-1},-4 - 4 t)\\
&(-4t^{-1},-4 t)&(-2t^{-1},-2 t)&(2t^{-1},2 t)&(4t^{-1},4 t)\\
&(-1 + s,-1 - s)&(1 + s,1 - s)&(-6 + 2 s,-6 - 2 s)&(-1 + 3 s,-1 - 3 s)\\
&(-4 + 4 s,-4 - 4 s)&&&\\
\hline \end{array} $

\end{center}
\label{data10b}
\end{table}

\begin{table}[htdp]
\caption{Boundary slope data for 2-bridge links with 10 crossings (part 3). The parameters $t$ and $s$ are rational with $1<t<\infty$ and $-1 \le s \le 1$.}
\begin{center}

$\begin{array}{|r|llll|}
\hline
p/q&\mbox{boundary slopes}&&&\\
\hline
\hline
23/64&(-2,-4)&(-2,-2)&(-2,0)&(0,-2)\\
&(0,0)&(4,0)&(4,2)&(-2-2t^{-1},-4 - 2 t)\\
&(4-2t^{-1},-2 - 2 t)&(-5-t^{-1},-3 - t)&(-5+t^{-1},-3 + t)&(-5+t^{-1},1 + t)\\
&(-2+2t^{-1},-2 + 2 t)&(-2+2t^{-1},2 + 2 t)&(4+2t^{-1},2 + 2 t)&(-4t^{-1},-2 - 4 t)\\
&(-2t^{-1},-2 - 2 t)&(-2t^{-1},-2 t)&(2t^{-1},2 t)&(2t^{-1},2 + 2 t)\\
&(4t^{-1},2 + 4 t)&(-5 + s,-5 - s)&(-3 + s,-3 - s)&(3 + s,3 - s)\\
&(5 + s,5 - s)&(2 + 2 s,2 - 2 s)&(-1 + 3 s,-1 - 3 s)&\\
25/66&(-4-t^{-1},-2 - t)&(8-t^{-1},-2 - t)&(8-t^{-1},2 - t)&(-4+t^{-1},-2 + t)\\
&(-4+t^{-1},6 + t)&(8+t^{-1},2 + t)&(-3t^{-1},-2 - 3 t)&(-t^{-1},-2 - t)\\
&(-t^{-1},-t)&(t^{-1},t)&(t^{-1},6 + t)&(3t^{-1},3 t)\\
&(3t^{-1},6 + 3 t)&(-4 + s,-4 - s)&(-2 + s,-2 - s)&(4 + 3 s,4 - 3 s)\\
&(6 + 3 s,6 - 3 s)&(2 + 5 s,2 - 5 s)&&\\
19/68&(-2,-2)&(-2,2)&(0,0)&(0,2)\\
&(2,0)&(2,4)&(4,2)&(-2-2t^{-1},-2 - 2 t)\\
&(4-2t^{-1},2 - 2 t)&(4-2t^{-1},-2 t)&(7-t^{-1},1 - t)&(7+t^{-1},3 + t)\\
&(-2+2t^{-1},2 t)&(-2+2t^{-1},4 + 2 t)&(4+2t^{-1},4 + 2 t)&(-4t^{-1},-4 t)\\
&(-2t^{-1},-2 t)&(2t^{-1},2 t)&(2t^{-1},2 + 2 t)&(4t^{-1},4 + 4 t)\\
&(2 s,-2 s)&(1 + s,1 - s)&(3 + s,3 - s)&(6 + 2 s,6 - 2 s)\\
&(3 + 3 s,3 - 3 s)&&&\\
29/70&(-11,-3)&(-6-t^{-1},-6 - t)&(-6-t^{-1},-2 - t)&(-6+t^{-1},-2 + t)\\
&(-3t^{-1},-8 - 3 t)&(-3t^{-1},-4 - 3 t)&(-3t^{-1},-3 t)&(-t^{-1},-4 - t)\\
&(-t^{-1},-t)&(t^{-1},t)&(3t^{-1},3 t)&(-5 + 2 s,-5 - 2 s)\\
&(-3 + 2 s,-3 - 2 s)&(-7 + 4 s,-7 - 4 s)&&\\
31/74&(-6-t^{-1},-2 - t)&(6-t^{-1},-4 - t)&(6-t^{-1},2 - t)&(-6+t^{-1},-2 + t)\\
&(-6+t^{-1},4 + t)&(6+t^{-1},2 + t)&(-3t^{-1},-4 - 3 t)&(-3t^{-1},-3 t)\\
&(-t^{-1},-4 - t)&(-t^{-1},-t)&(t^{-1},t)&(t^{-1},4 + t)\\
&(3t^{-1},3 t)&(3t^{-1},4 + 3 t)&(5 s,-5 s)&(-5 + 2 s,-5 - 2 s)\\
&(-3 + 2 s,-3 - 2 s)&(3 + 2 s,3 - 2 s)&(5 + 2 s,5 - 2 s)&\\
21/76&(-6,-4)&(-6,2)&(-2,-2)&(-2,2)\\
&(0,-2)&(0,0)&(4,-2)&(4,0)\\
&(4,2)&(-4-2t^{-1},-2 - 2 t)&(-2-2t^{-1},-4 - 2 t)&(4-2t^{-1},-4 - 2 t)\\
&(4-2t^{-1},-2 t)&(-4+2t^{-1},-2 + 2 t)&(-4+2t^{-1},2 + 2 t)&(-2+2t^{-1},-2 + 2 t)\\
&(-2+2t^{-1},2 + 2 t)&(4+2t^{-1},2 + 2 t)&(-4t^{-1},-2 - 4 t)&(-2t^{-1},-2 t)\\
&(2t^{-1},2 t)&(2t^{-1},2 + 2 t)&(4t^{-1},2 + 4 t)&(-5 + s,-5 - s)\\
&(-1 + s,-1 - s)&(2 s,-2 s)&(3 + s,3 - s)&(5 + s,5 - s)\\
&(2 + 2 s,2 - 2 s)&(1 + 3 s,1 - 3 s)&(-2 + 4 s,-2 - 4 s)&\\
31/80&(-8,-4)&(-8,0)&(-4,-2)&(-4,0)\\
&(0,0)&(2,-2)&(2,2)&(-4-2t^{-1},-4 - 2 t)\\
&(2-2t^{-1},-6 - 2 t)&(2-2t^{-1},-2 - 2 t)&(2-2t^{-1},2 - 2 t)&(-4+2t^{-1},2 t)\\
&(2+2t^{-1},2 + 2 t)&(-4t^{-1},-4 - 4 t)&(-2t^{-1},-2 - 2 t)&(2t^{-1},2 t)\\
&(4t^{-1},4 t)&(-3 + s,-3 - s)&(2 s,-2 s)&(-6 + 2 s,-6 - 2 s)\\
&(-2 + 2 s,-2 - 2 s)&(-4 + 4 s,-4 - 4 s)&&\\
\hline \end{array} $

\end{center}
\label{data10c}
\end{table}

 \end{document}